\documentclass[letterpaper,10pt]{article}
\pdfoutput=1
\usepackage{mathtools}
\usepackage{amssymb}
\usepackage{mathrsfs,dsfont,textcomp}

\usepackage{empheq}

\usepackage{amsthm}
\usepackage{thmtools} 

\usepackage{titlesec}  

\usepackage{parskip}
		
\usepackage{enumitem}

\usepackage{eso-pic}
\usepackage{tikz}
\usetikzlibrary{tikzmark}
\usetikzlibrary{angles, quotes}
\usetikzlibrary{intersections}
\usetikzlibrary{calc,decorations.markings,arrows.meta}

\usepackage{footnote}
\usepackage{footmisc} 
\newcounter{savefootnote}
\newcounter{symfootnote}
\newcommand{\symfootnote}[1]{
   \setcounter{savefootnote}{\value{footnote}}
   \setcounter{footnote}{\value{symfootnote}}
   \ifnum\value{footnote}>8\setcounter{footnote}{0}\fi
   \let\oldthefootnote=\thefootnote
   \renewcommand{\thefootnote}{\fnsymbol{footnote}}
   \footnote{#1}
   \let\thefootnote=\oldthefootnote
   \setcounter{symfootnote}{\value{footnote}}
   \setcounter{footnote}{\value{savefootnote}}
}

\usepackage{etoolbox}
\patchcmd{\thebibliography}{\section*{\refname}}{}{}{}

\usepackage[letterpaper, top=3cm, bottom=2cm, left=2.5cm, right=2.5cm]{geometry}

\titleformat
{\section}
{\scshape} 
{\thesection.}  
{1.0pc} 
{} 
\titlespacing  
{\section}
{0pc}  
{1.5pc} 
{\bigskipamount} 

\titleformat
{\subsection}
{\bfseries}
{\thesubsection}
{1.0pc}
{}
\titlespacing
{\subsection}
{0pc}
{\bigskipamount}
{\medskipamount}

\declaretheoremstyle[
spaceabove=\medskipamount, spacebelow=\medskipamount,
headfont=\normalfont\bfseries,
bodyfont=\normalfont,
postheadspace=  12pt,]
{mystyledefinition}
\declaretheorem[style=mystyledefinition, numberwithin=section]{definition}

\declaretheoremstyle[
spaceabove=\medskipamount, spacebelow=\medskipamount,
headfont=\normalfont\sl,
bodyfont=\normalfont,
postheadspace=  12pt,]
{mystyleexample}
\declaretheorem[style=mystyleexample, numberwithin=section, numberlike=definition]{example}
\declaretheorem[style=mystyleexample, numberwithin=section, numberlike=definition]{remark}

\declaretheoremstyle[
spaceabove=\medskipamount, spacebelow=\medskipamount,
headfont=\normalfont\bfseries,
bodyfont=\sl,
postheadspace=  12pt,  
]{mystylethm}
\declaretheorem[style=mystylethm, numberwithin=section, numberlike=definition]{theorem}
\declaretheorem[style=mystylethm, numberwithin=section, numberlike=definition]{proposition}
\declaretheorem[style=mystylethm, numberwithin=section, numberlike=definition]{lemma}
\declaretheorem[style=mystylethm, numberwithin=section, numberlike=definition]{corollary}

\font\titlefont= cmbx12
\font\authorfont=cmcsc10 
\font\headfont=cmcsc10

\font\small=cmr9
\font\smallbf=cmbx9

\font\smalltt=cmtt9

\def\Proof{\noindent {\sl Proof.}\enspace}

\def\qedmark{\hbox{\vrule height 4pt width 3pt}}
\def\qedskip{\vrule height 4pt width 0pt depth 1pc}
\def\qed{\penalty 1000\quad\penalty 1000{\qedmark\qedskip}}

\def\bx{\hbox{\bf x}}

\def\bz{\hbox{\bf z}}
\def\bv{\hbox{\bf v}}

\def\bT{\hbox{\bf T}}
\def\bS{\hbox{\bf S}}

\def\bK{\hbox{\bf K}}

\def\cA{{\cal A}}
\def\cC{{\cal C}}
\def\cD{{\cal D}}
\def\cF{{\cal F}}
\def\cE{{\cal E}}
\def\cI{{\cal I}}
\def\cJ{{\cal J}}
\def\cK{{\cal K}}
\def\cN{{\cal N}}
\def\cO{{\cal O}}
\def\cQ{{\cal Q}}

\def\cV{{\cal V}}

\def\srQ{\mathscr{Q}}

\def\gF{\mathfrak{F}}
\def\gI{\mathfrak{I}}

\def\bxi{\boldsymbol{\xi}}

\def\R{{\rm I \! R}}

\renewcommand{\r}{\rightarrow}

\begin{document}

\begin{center}
{\titlefont Maslov Index, Spectral Flow and Bifurcation of Electromagnetic \par Geodesics Within an Energy Level}
\end{center}
\bigskip
\centerline{{\authorfont Henrique Vit\'orio}\symfootnote{Departamento de Matem\'atica, 
Universidade Federal de Pernambuco, Brazil. E-mail address: \smalltt henrique.vitori@ufpe.br} }

\vskip .5cm

\bigskip
\medskip
\hfil {\hsize = 10cm \vbox{\noindent {\smallbf Abstract.} \small
We develop appropriate notions of Maslov index and spectral flow for electromagnetic geodesics within a fixed energy level
and prove a Morse Index type theorem in this context. This is then applied to the  
problem of electromagnetic geodesics, all of whose energy are the same, bifurcating from
a given one.}} \hfil

\vskip.3cm

\bigskip

\noindent 

\section{Introduction}

\noindent Let $M$ be a $n-$dimensional smooth manifold endowed with both a semi-Riemannian metric $g$ and a closed
$2-$form $\sigma$. Let also $Y:TM\r TM$ be the endomorphism field on $M$ defined via 
$\sigma(p)[u,v]=g(p)[u,Y(p)[v]]$. We consider in this work the following system of second order differential equations on $M$,
\begin{equation}\label{lorentzforceequation}
\frac{D}{dt} \dot{\gamma} = Y(\gamma)[\dot{\gamma}],
\end{equation}
for $D/dt$ the semi-Riemannian covariant derivative. In case $g$ is Riemannian, resp. Lorentzian, these equations govern
the motion of a particle on $(M,g)$ under 
the influence of the magnetic, resp. electromagnetic, field $\sigma$. In the general case of an arbitrary semi-Riemannian 
metric $g$, we shall use 
the term {\sl electromagnetic}. The endomorphism $Y$ is known as the Lorentz force of the electromagnetic field 
$\sigma$ and
Eq.\ (\ref{lorentzforceequation}) as the {\sl Lorentz force equation}.  

We shall call solutions to Eq.\ (\ref{lorentzforceequation}) {\sl electromagnetic geodesics}. As it is well-known, along an electromagnetic
geodesic $\gamma:I\subseteq\R \r M$ the {\sl energy} 
\[ 
E(\gamma,\dot{\gamma})=\frac{1}{2}g(\gamma)[\dot{\gamma},\dot{\gamma}] 
\]
stays constant. In this work we are interested in the following phenomenon:  

\begin{definition}\label{definitionbifurcation}
Let $\gamma:[0,T] \r M$ be an electromagnetic geodesic with energy $\kappa$. We say that an instant 
$t_0\in (0,T)$ is an {\sl energy-constrained bifurcation instant} for $\gamma$ if there exists a sequence 
$\gamma_n:[0, T]\rightarrow M$ of distinct electromagnetic geodesics
and two sequences of instants $(t_n)_{n\geq 1}$ and $(t_n')_{n\geq 1}$ in $(0,T)$ 
such that
\begin{enumerate}[left=4pt .. \parindent]
\item $\gamma_n(0)=\gamma(0)$ for all $n$, and $\dot{\gamma}_n(0) \r \dot{\gamma}(0)$ as 
$n \r \infty$;
\item $t_n \r t_0$ and $t_n' \r t_0$ as $n\r \infty$;
\item $\gamma_n(t_n')=\gamma(t_n)$ for all $n$; 
\item $\gamma_n$ has energy $\kappa$ for all $n$.
\end{enumerate}
\end{definition}

This kind of qualitative behavior of solutions to variational problems arising in differential geometry, 
along with its variations, has been the object of study 
of many works, including \cite{tausk2004,piccione_nonlinearity,portaluri_piccione,portaluri2007,giambo}.
Bifurcation phenomena are natural to consider. For example, a slight variation of the above definition, in the 
context of light-like geodesics of a Lorentzian manifold, can be used to model the relativistic phenomenon of
{\sl gravitational lensing} (see \cite{piccione_nonlinearity}). Considering a bifurcation perspective has also provided 
valuable insights to the study of solutions to some important problems in multidimensional variational calculus; for
an updated overview, we refer to the survey \cite{piccione_notices}.  

The novelty in our work is the energy constraint expressed in 4.\ of Definition \ref{definitionbifurcation}. In contrast, the
definition of bifurcation in \cite{portaluri_piccione}, which
we shall refer here as {\sl ordinary} bifurcation, while does not require the energy constraint, it demands
that the two sequences $(t_n)$ and $(t_n')$ in Definition \ref{definitionbifurcation} can be made the same.

\begin{example}\label{examplebifurcation1}
Consider a constant magnetic field $\sigma$ on the flat $\R^2$ given by the area form. According to our sign convention,
the Lorentz force $Y$ is given by multiplication by $-i$ (via the identification $\R^2=\mathds{C}$), so that the Lorentz
force equation becomes $\ddot{\gamma}=-i\dot{\gamma}$. The solutions to this equation are the circles parametrized clockwise 
with unit angular velocity, and the energy of such a circle is thus half the square of its radius. 
Given such a circle $\gamma:[0,2\pi]\r \R^2$, by letting it ``hang'' from the point $\gamma(0)$, without changing its radius,
one easily sees that $t_0=\pi$ is an 
energy-constrained bifurcation instant along $\gamma$; see the figure below. 
\begin{figure}[!h]
\centering
\begin{tikzpicture}[scale=0.5,
    upper clockwise arrow/.style={
        postaction={decorate},
        decoration={
            markings,
            mark=at position 0.25 with {\arrowreversed[thin]{Stealth[scale=0.8]}}
        }
    },
    upper anticlockwise arrow/.style={
        postaction={decorate},
        decoration={
            markings,
            mark=at position 0.25 with {\arrow[thin]{Stealth[scale=0.8]}}
        }
    }
]

\newcommand{\radius}{4}

\draw[thick, black, upper clockwise arrow, name path=circleC] (0,0) circle (\radius);

\draw[thin, red, upper clockwise arrow, name path=circleC1] (-0.45, 1.85) circle (\radius);

\draw[thin, red, upper clockwise arrow, name path=circleC2] (-0.24,1.36) circle (\radius);

\draw[thin, red, upper clockwise arrow, name path=circleC3] (-0.13,1) circle (\radius);

\draw[thin, red, upper clockwise arrow, name path=circleC4] (-0.05,0.61) circle (\radius);

\fill[black, name intersections={of=circleC and circleC1, name=i1, total=\t}]
    \foreach \s in {1,...,\t} {(i1-\s) circle (2pt)};
  
\fill[black, name intersections={of=circleC and circleC2, name=i2, total=\t}]
    \foreach \s in {1,...,\t} {(i2-\s) circle (2pt)};

\fill[black, name intersections={of=circleC and circleC3, name=i3, total=\t}]
    \foreach \s in {1,...,\t} {(i3-\s) circle (2pt)};

\fill[black, name intersections={of=circleC and circleC4, name=i4, total=\t}]
   \foreach \s in {1,...,\t} {(i4-\s) circle (2pt)};
   
\node[below right, yshift=0pt, xshift=3pt] at (i2-1) {$\gamma_n(t_n')=\gamma(t_n)$};

\fill[black] (4,0) circle (2pt);
\node[below, yshift=-2pt, right] at (4,0) {$\gamma(\pi)$};

\fill[black] (-4,0) circle (2pt);
\node[below, yshift=-2pt, left] at (-4,0) {$\gamma(0)$};

\node[below] at (0,4) {$\gamma$};

\node[below, red] at (0,4.8) {$\gamma_n$};  
\end{tikzpicture}
\end{figure}
Also, even by allowing the radius of $\gamma_n$ to vary, the orbits $\gamma_n$ and $\gamma$ will 
take different times to reach their intersection because they have the same period. 
Thus, $t_0=\pi$ is not
an ordinary bifurcation instant along $\gamma$. On the other hand, it is clear that $t_0=2\pi$ is both energy-constrained and 
ordinary bifurcation instant
along $\gamma$ (in this case, necessarily $t_n'=t_n=2\pi$ for all $n$). 

Let us now move the setup to the round sphere $S^2$ with radius 1 endowed
with a uniform magnetic field $\sigma = \beta dS^2$, where $dS^2$ is the area form and $\beta\neq 0$ is a constant.
The Lorentz force equation then writes 
\begin{equation}\label{eqlorentzsphere}
\frac{D}{dt}\dot{\gamma} = -\beta (\gamma \times \dot{\gamma}).
\end{equation} 
The solutions to this equation are the circles on $S^2$ parametrized with constant angular velocity whose 
radius $r\in (0,1)$ and angular velocity $\omega$ 
are tied by the relation $\omega\sqrt{1-r^2} = -\beta$. Thus their period decreases with the energy (hence with the radius). 
This ensures that nearby orbits arrive at the same time at their intersection (up to first order) giving rise to an ordinary bifurcation
instant that is not energy-constrained; we shall do this analytically in Example \ref{exampleconjugate}. 
\end{example}
  
The energy constraint is a very natural one since
the dynamics of a Hamiltonian flow can change significantly from one energy level to another.
Moreover, such a constraint may convey a physical meaning
to the problem. For example, for $(M,g)$ a space-time,
the constraint $E(\gamma_n,\dot{\gamma}_n)=-1/2$ means that the trajectories $\gamma_n$ are parametrized by
{\sl proper time}; 
we remark that, in this context, a bifurcation result for timelike electromagnetic geodesics parametrized by proper time, and departing
from a timelike curve representing an observer, was obtained in \cite{giambo} by exploring the 
{\sl Kaluza-Klein} correspondence and resorting to (a stronger form of) the bifurcation
result for lightlike geodesics in \cite{piccione_nonlinearity}. 

In the absence of the energy constraint,  
sufficient conditions for 
existence of bifurcation have been obtained in \cite{tausk2004}, \cite{portaluri_piccione} and \cite{portaluri2007} for the 
Euler-Lagrange equations associated to the {\sl geodesic}, {\sl electromagnetic} and {\sl perturbed geodesic} Lagrangians $L:TM\r \R$, respectively. 
These were achieved
by proving a Morse Index type theorem and then applying to the action functional of $L$
the functional analytic 
results for bifurcation of critical points of strongly indefinite functionals established in \cite{recht}.
Observe that since we are allowing for a {\sl non-exact}  electromagnetic field $\sigma$, 
Eq.\ (\ref{lorentzforceequation}) might not be the Euler-Lagrange equations of some Lagrangian function 
on $TM$. Nevertheless, the analysis of the electromagnetic geodesics in the vicinity of a given one
$\gamma:[0,T]\r M$ connecting two points $p$ and $q$, and having the same energy as $\gamma$, can also be brought to the realm of critical point theory thanks
to the following (see \cite{gabrielle} and Sec.\ \ref{sectionvariationalsetup} below): {\sl the solutions $\gamma:[0,T]\r M$ to Eq.\ (\ref{lorentzforceequation}) 
connecting two points $p$ and $q$ and constrained to $E(\gamma,\dot{\gamma})=\kappa$ correspond to the zeros of a 
certain smooth $1-$form $\eta_\kappa$ on $\Omega_{p,q}([0,1])\times \R_+$, where
$\Omega_{p,q}([0,1])$ is the Hilbert manifold of (Sobolev regular) curves $x:[0,1]\r M$ connecting $p$ and $q$. Furthermore, 
$\eta_\kappa$ possesses local primitives.}

Our purpose is to provide a general sufficient condition for energy-constrained bifurcation to occur by establishing 
a Morse Index type theorem in the energy-constrained setting.
Each electromagnetic geodesic $\gamma:[0,T]\r M$ with non null energy $\kappa$ possesses an {\sl energy-constrained index form}, which
is the quadratic form $\cQ_{\kappa,\sigma}(\gamma)$ given by the second derivative of a local primitive for the $1-$form $\eta_\kappa$ 
at the critical point corresponding to $\gamma$. 
By analysing the kernel of $\cQ_{\kappa,\sigma}(\gamma)$ we are led to the following notion
of conjugacy along $\gamma$ (we remark that the same notion of conjugacy was obtained in \cite{giambo} by different means); 
in the following, the symbols $\nabla$ and $R$ will denote the Levi-Civita 
connection and the curvature tensor of the metric $g$, respectively, with the sign convention $R(X,Y)=[\nabla_X,\nabla_Y]-\nabla_{[X,Y]}$.

\begin{definition}
Let $\gamma:[0, T]\r M$ be an electromagnetic geodesic with energy $\kappa\neq 0$. An {\sl energy-constrained Jacobi field}
along $\gamma$ is a (smooth) vector field $J$ along $\gamma$ that satisfies the equation
\begin{equation}\label{eqjacobiintro}
\frac{D^2}{dt^2} J + R(J,\dot{\gamma})\dot{\gamma} - (\nabla_J Y)[\dot{\gamma}] - Y \left[ \frac{D}{dt}J \right]
-\frac{1}{2\kappa} g\left[ \frac{D}{dt}J(0) , \dot{\gamma}(0) \right] Y[\dot{\gamma}] = 0.
\end{equation}
An instant $t_0\in (0,T]$, for which there exists a non-null energy-constrained Jacobi field that
vanishes at both $t=0$ and $t=t_0$ will be called an {\sl energy-constrained conjugate instant} along $\gamma$.  
\end{definition} 

With a different terminology, Eq.\ (\ref{eqjacobiintro}) and the above notion of conjugacy first appeared in \cite{giambo} where
they played a major role in establishing the aforementioned authors' bifurcation result.

Without the energy constraint, we shall use the term {\sl ordinary} to refer to the well-known analogous concepts. So, for instance, we shall
also speak of ordinary Jacobi fields and ordinary conjugacy; these are obtained by erasing the last term of the left-hand side of Eq.\ (\ref{eqjacobiintro}).  
It turns out that our definition of conjugacy is equivalent to the one presented in \cite{assenza2} (see also \cite{reber_shen}) in terms 
of ordinary Jacobi fields: in that work, $t_0$ is said a conjugate instant along $\gamma$ if there exists
a non null ordinary Jacobi field $J$ satisfying $g(\gamma)\left[DJ/dt,\dot{\gamma}\right]\equiv 0$  
and that vanishes at $t=0$ and is tangent to $\gamma$ at $t=t_0$; the equivalence of both definitions is shown in Sec.\ \ref{subsectionversus}.

Let us suppose for the rest of this introduction that $\gamma:[0,T]\r M$ is an electromagnetic geodesic with
{\sl non null} energy $\kappa$.
An application of the Implicit Function
Theorem (see the analogous remark in Sec.\ \ref{sectionabstractbifurcation}), together with a bundle trivialization argument carried out
in Sec.\ \ref{subsectionproofbifurcation}, shows that every energy-constrained bifurcation instant is an energy-constrained conjugate instant. 
Reciprocally, Theorem \ref{theorembifurcationmaslov} below establishes that the interval $(0,T]$ will contain an 
energy-constrained bifurcation instant for $\gamma$ provided that
the ``net number'' of energy-constrained conjugate instants in $(0,T]$ be non null. This net number shall be given by the following
notion of {\sl energy-constrained Maslov index} of $\gamma$. For this, let  $\zeta_t : TM \rightarrow TM$ be the flow corresponding to Eq.\ (\ref{lorentzforceequation}). 
This flow is the Hamiltonian flow of the energy function $E:TM\rightarrow \R$ when one endows $TM$ with the  
symplectic structure $\omega$ induced 
by both the metric $g$ and the electromagnetic field $\sigma$, i.e. $\omega = d \alpha_g + \pi^* \sigma$,
for $\alpha_g$ the pull-back of the canonical $1-$form $\alpha$ of $T^*M$ by the diffeomorphism 
$g^\sharp:TM\r T^*M$ and $\pi:TM\r M$ the projection map.
Let $S$ be the corresponding Hamiltonian vector field. Since 
the energy surface $E^{-1}(\kappa)$ is invariant by $\zeta_t$,  
the derivative of $\zeta_t$ induces an isomorphism of quotient spaces, 
\begin{equation}\label{mapbarphi}
\bar{\Phi}_t: T_{(p,v)} E^{-1}(\kappa)/\left\langle S \right\rangle \, 
\rightarrow \,
T_{(\gamma,\dot{\gamma})} E^{-1}(\kappa)/\left\langle S \right\rangle ,
\end{equation}
where $p=\gamma(0)$ and $v=\dot{\gamma}(0)$. The above
quotient spaces inherit from $\omega$ a symplectic structure, and the map $\bar{\Phi}_t$ is symplectic. Let 
$V E^{-1}(\kappa)\subset T E^{-1}(\kappa)$ denote the vertical distribution on $E^{-1}(\kappa)$, which is given by the tangent
spaces to the fibers of the projection $E^{-1}(\kappa) \rightarrow M$. The image of $V E^{-1}(\kappa)$ in 
$T E^{-1}(\kappa)/\langle S \rangle$ is a Lagrangian subbundle, which we shall denote by $\bar{L}$. 
By pulling $\bar{L}_{(\gamma,\dot{\gamma})}$ back by $\bar{\Phi}_t$ we thus obtain a curve of Lagrangian subspaces,
\begin{equation}\label{lagrangianpath}
\ell: [0,T] \rightarrow \Lambda\left( T_{(p,v)}E^{-1}(\kappa)/\langle S \rangle \right), \quad 
\ell(t) = (\bar{\Phi}_t)^{-1} \bar{L}_{(\gamma,\dot{\gamma})}.
\end{equation}
In the above, $\Lambda(-)$ denotes the Lagrangian Grassmannian manifold of all Lagrangian subspaces
of the corresponding symplectic vector space. 
In the same way as ordinary conjugate instants correspond 
to non null intersections of the pull-back $\check{\ell}(t)$ of the vertical distribution $VTM$ 
with the fixed vertical space $V_{(p,v)}TM$,
it turns out that energy-constrained conjugate instants correspond to non null intersections of $\ell(t)$ with $\bar{L}_{(p,v)}$.
Indeed, we show in Sec.\ \ref{sectionsymplectic} that, up to a symplectic change of variables, the ``flow'' of the energy-constrained Jacobi equation
(\ref{eqjacobiintro}) is nothing but $\bar{\Phi}_t$. 
We remark that this correspondence between energy-constrained conjugate instants and non null intersections $\ell(t)\cap \bar{L}_{(p,v)}$ has
also been noted in \cite{assenza2} (see Proposition 4.3 of \cite{assenza2}) where it played a role
in the authors' construction of the so-called {\sl magnetic Green bundles}. It is interesting to note that 
the definition of conjugacy in \cite{assenza2} directly amounts to requiring that $\check{\ell}(t)\cap L_{(p,v)}\neq\{0\}$,
where $L$ is the Lagrangian distribution along $E^{-1}(\kappa)$ given by $L=\langle S \rangle + VTM \cap T E^{-1}(\kappa)$;
on the other hand, by a symplectic reduction argument, one has $\check{\ell}(t)\cap L_{(p,v)}\neq\{0\}$ if, and only if, $\ell(t)\cap \bar{L}_{(p,v)}\neq \{0\}$.

\begin{definition}
The {\sl energy-constrained Maslov index} of $\gamma$, denoted by $\mu_\kappa(\gamma)$, is defined as the 
Maslov index of the Lagrangian path $\ell|_{[\epsilon, T]}$ with respect to the reference Lagrangian 
$\bar{L}_{(p,v)}$, for $\epsilon$ any sufficiently small positive number (in the sense of Proposition \ref{propisolatedinstant}). 
\end{definition} 

On the other hand, by traveling along $\gamma$ from $p=\gamma(0)$ to $q=\gamma(T)$, the resulting family of energy-constrained index forms 
can be identified with a path of quadratic forms $Q=\{Q_s\}_{s\in (0,T]}$ on the Hilbert space
$H_0^1([0,1],\R^n)\times \R$ which are of Fredholm type. For such a family one can compute the so-called  
{\sl spectral flow} along an interval $[a,b]\subset (0,T]$, denoted by ${\rm sf}(Q,[a,b])$. 
Roughly, it gives the net number of the eigenvalues
of the self-adjoint operator $L_s$ representing $Q_s$ that cross zero as $s$ runs from $a$ to $b$. 
The spectral flow is an important homotopy invariant which substitutes the notion of Morse index when the latter fails to
be finite, and plays a key role in strongly-indefinite problems. For instance,  
the spectral flow is a key to understanding variational bifurcation of strongly-indefinite functionals. Indeed, the main result of
\cite{recht} shows that the non vanishing of the spectral flow implies bifurcation. 
We define the {\sl energy-constrained spectral flow} of $\gamma$, and denote it by ${\rm sf}_\kappa(\gamma)$, as the spectral flow
${\rm sf}(Q,[\epsilon,T])$ for any sufficiently small positive $\epsilon$ (as in Proposition \ref{propisolatedinstant}).   
The Morse Index type theorem that we shall prove establishes that the spectral flow, which is infinite-dimensional in nature, can 
be computed as the more tractable Maslov index:

\begin{theorem}\label{morseindextheorem}
Suppose that $T$ is neither an energy-constrained nor an ordinary conjugate instant along $\gamma$. Then, we have an equality
\begin{equation}\nonumber
{\rm sf}_\kappa(\gamma) = -\mu_\kappa(\gamma).
\end{equation}
\end{theorem}  
This theorem is the energy-constrained version of the Morse Index theorem of \cite{portaluri_piccione} which establishes an 
analogous equality for the {\sl ordinary} Maslov index $\mu(\gamma)$ and spectral flow ${\rm sf}(\gamma)$ of an electromagnetic geodesic. 
Its proof shall consist of three steps. First we compute the difference ${\rm sf}_\kappa(\gamma)-{\rm sf}(\gamma)$ by establishing 
in Sec.\ \ref{sectionrestrictionspectralflow} a general formula relating the spectral flow of a path $\{\cQ_s\}$ to the spectral flow of
the path obtained by restricting the $\cQ_s$ to a fixed closed finite codimensional subspace of the ambient Hilbert space. This result extends 
to general paths (but with invertible endpoints), and with a considerably simpler proof, Theorem 4.4 of \cite{benevieri} proved for paths
that are compact perturbations of a {\sl fixed} symmetry, and might be of independent interest; another related formula has also been 
obtained in \cite{portaluri_wu_yang}. 
It is important to note that the difference ${\rm sf}_\kappa(\gamma)-{\rm sf}(\gamma)$, which is either $0$ or $\pm 1$, has been the subject of 
investigation of some works dealing with index computations of periodic orbits, in which a formula for such difference was expressed in terms 
of the concept of {\sl orbit cylinder}: in \cite{paternain_merry} and \cite{paternain2015} for the case
of Tonelli Lagrangian systems, in which case the spectral flows reduce (up to sign) to the Morse indices of the free 
period and fixed period action functionals; and in \cite{portaluri_wu_yang} for non-convex Lagrangian systems.   
For the sake of comparison, we show in Sec.\ \ref{subsectionorbitcylinder} how our computation expresses in terms
of a corresponding notion of orbit cylinder, called here {\sl orbit strip}. In the second step of the proof, by regarding the Lagrangian
path $\ell(t)$ as the symplectic reduction of the Lagrangian path $\check{\ell}(t)$ on $\Lambda(T_{(p,v)}TM)$,
we resort to a result from \cite{vitorio} to equate $\mu_\kappa(\gamma)$ to a Maslov index of $\check{\ell}(t)$. 
This has the advantage that the latter Maslov index differs
from $\mu(\gamma)$ by a {\sl Hörmander four-fold index}, which in turn is half the difference of two {\sl Kashiwara indices}
of triplets of Lagrangian subspaces. So step three shall consist of computing the Kashiwara indices of the triplets formed 
by the endpoints of the path $\check{\ell}(t)$ and the reference Lagrangians. With the computations of ${\rm sf}_\kappa(\gamma)-{\rm sf}(\gamma)$ and
$\mu_\kappa(\gamma)-\mu(\gamma)$ at hand, Theorem \ref{morseindextheorem} will follow from the Morse index theorem of \cite{portaluri_piccione}.
We remark that works 
\cite{paternain2015} and \cite{portaluri_wu_yang} also relate the Morse indices/spectral flows with notions of Maslov index, although their 
definition of the latter does not take into account the energy constraint.     

As a consequence of Theorem \ref{morseindextheorem} and the criterion for variational bifurcation of 
strongly-indefinite functionals from \cite{recht}, we shall obtain the following criterion for energy-constrained 
bifurcation of electromagnetic geodesics.  
  
\begin{theorem}\label{theorembifurcationmaslov}
Suppose that $T$ is neither an energy-constrained nor an ordinary conjugate instant along $\gamma$. If 
$\mu_\kappa(\gamma)\neq 0$, then there exists an energy-constrained bifurcation instant for $\gamma$ in $(0,T]$. 
\end{theorem}

For an isolated energy-constrained conjugate instant $t_0$ along $\gamma$, the energy-constrained Maslov index of $\gamma$ {\sl across} $t_0$, 
denoted by $\mu_\kappa(\gamma,t_0)$, is defined as the Maslov index of the Lagrangian path $\ell|_{[t_0-\epsilon , t_0+\epsilon]} $ with respect to the Lagrangian  
$\bar{L}_{(p,v)}$, for $\epsilon$ a sufficiently small positive number. 

\begin{corollary}\label{corollarymaslovacross}
An isolated energy-constrained conjugate instant $t_0$ will be an energy-constrained bifurcation instant provided that $\mu_\kappa(\gamma,t_0)\neq 0$.
\end{corollary}

For an energy-constrained conjugate instant $t_0$ we consider the subspace $\mathds{J}'[t_0]$ of $T_{\gamma(t_0)}M$ formed by the
values $DJ/dt|_{t=t_0}$, for $J$ varying among the energy-constrained Jacobi fields vanishing at both $t=0$ and $t=t_0$. Also, let 
$\mathfrak{g}_{t_0}$ be the quadratic form on $T_{\gamma(t_0)}M$ given by
\begin{equation}\nonumber
\mathfrak{g}_{t_0}[u,v] = g(\gamma(t_0))[u,v] - \frac{1}{2\kappa} g(\gamma(t_0))\bigl[\dot{\gamma}(t_0) , u \bigr] g(\gamma(t_0))\bigl[ \dot{\gamma}(t_0) , v \bigr].
\end{equation}  
We say that $t_0$ is {\sl non-degenerate} if the restriction of $\mathfrak{g}_{t_0}$ to $\mathds{J}'[t_0]$ is non-degenerate.
For such an instant we shall show in Proposition \ref{propsignaturemaslov} that $\mu_\kappa(\gamma,t_0)$ is equal to 
the signature (see (\ref{signaturedefinition}) below) of the restriction of $\mathfrak{g}_{t_0}$ to $\mathds{J}'[t_0]$. Therefore,  

\begin{corollary}\label{corollarybifurcation}
A non-degenerate energy-constrained conjugate instant along $\gamma$ with ${\rm sig}\left(\mathfrak{g}_{t_0}|_{\mathds{J}'[t_0]}\right)\neq 0$ 
is an energy-constrained bifurcation instant. 
\end{corollary}

The paper is organized as follows: in Sec.\ \ref{sectionmaslov} we recall the definition and basic properties of the Maslov index. In particular, 
in Sec.\ \ref{sectionkashiwara} we define the Kashiwara index and sketch the proof of its relation with the difference of Maslov indices. 
The general notion of spectral flow, as well as of variational bifurcation, is discussed in Sec.\ \ref{sectionspectralflow}. In particular, we prove in
Sec.\ \ref{sectionrestrictionspectralflow} a result for the spectral flow of a restriction to a closed finite codimensional subspace. 
The necessary variational setup for the problem of electromagnetic geodesics is established in Sec.\ \ref{sectionvariationalsetup}, and
energy-constrained Jacobi fields and conjugacy are 
considered in Sec.\ \ref{sectionjacobifields}. In Sec.\ \ref{sectionsymplectic} we describe the energy-constrained conjugacy under a symplectic point of
view by first identifying a natural pre-symplectic structure which is preserved by the energy-constrained Jacobi equation, and then 
showing that a quotient of that equation is symplectically equivalent to the flow (\ref{mapbarphi}). The definition of the energy-constrained
spectral flow of an electromagnetic geodesic is carried out in Sec.\ \ref{sectionanalytic} along with the proof of the necessary
analytic facts. At last, the proofs of Theorem \ref{morseindextheorem}, Theorem \ref{theorembifurcationmaslov} and 
Corollary \ref{corollarymaslovacross} are carried out in Sec.\ \ref{sectionproofs}.

\section{The Maslov index of a curve of Lagrangian subspaces}\label{sectionmaslov}

\noindent Let $(V,\omega)$ be a finite-dimensional real symplectic vector space. The manifold of Lagrangian
subspaces of $V$, the {\sl Lagrangian Grassmannian}, shall be denoted by $\Lambda (V)$. We shall here recall
the definition of the Maslov index of a continuous path $\ell:[a,b]\rightarrow \Lambda(V)$ with respect to
a given reference Lagrangian $L_0\in\Lambda(V)$ through a uniqueness result from \cite{maslovbook}.

Given a pair $(L_0,L_1)$ of transverse Lagrangian subspaces, a chart for $\Lambda(L)$ around $L_0$ is obtained as follows. 
Firstly, given $L\in\Lambda(V)$, let us
denote $\Lambda^0(L)=\left\{ L'\in\Lambda(V) \; : \; L\cap L'=\{0\} \right\}$; this is an open and dense subset
of $\Lambda(V)$. Now, every $L\in \Lambda^0(L_1)$ is the graph of a unique linear map
$T:L_0 \rightarrow L_1$ such that the bilinear form on $L_0$ defined by $\omega\left[ T\cdot  ,  \cdot  \right]|_{L_0}$ 
is symmetric. A chart around $L_0$ is then obtained by setting
\begin{equation}\nonumber
\left\{
\begin{aligned}
& \varphi_{L_0,L_1}: \Lambda^0(L_1) \rightarrow {\rm B}_{\rm sym}(L_0), \\
& \varphi_{L_0,L_1}(L)= \omega\left[ T\cdot  , \cdot  \right]|_{L_0},
\end{aligned}\right.
\end{equation}
where ${\rm B}_{\rm sym}(L_0)$ is the space of symmetric bilinear forms on $L_0$. The Maslov index 
of $\ell:[a,b]\rightarrow \Lambda(V)$ with respect to $L_0$ is the semi-integer $\mu_{L_0}(\ell)$ 
characterized by the following uniqueness result (Corollary 5.2.4 of \cite{maslovbook}). Before proceeding, recall
that the {\sl signature} of a (real) quadratic form $Q$, denoted ${\rm sig}\,Q$, is the difference between
the indices of $-Q$ and $Q$,
\begin{equation}\label{signaturedefinition}
{\rm sig}\, Q = {\rm ind}(-Q) - {\rm ind}\, Q.
\end{equation}

\begin{theorem}\label{theoremmaslovindexdefinition}
Given a Lagrangian $L_0\in\Lambda(V)$, there exists a unique function that assigns to each continuous path
$\ell:[a,b]\rightarrow\Lambda(V)$ a semi-integer $\mu_{L_0}(\ell)$ and that satisfies 
\begin{enumerate}[left=4pt .. \parindent]
\item $\mu_{L_0}$ is invariant under homotopies that fix endpoints,
\item $\mu_{L_0}$ is additive by concatenation,
\item if $\ell([a,b])\subset \Lambda^0(L_1)$ for some $L_1\in \Lambda^0(L_0)$, then
\begin{equation}\label{eq2931}
\mu_{L_0}(\ell) = \frac{1}{2}{\rm sig}\, \varphi_{L_0,L_1}(\ell(b)) - \frac{1}{2}{\rm sig}\,\varphi_{L_0,L_1}(\ell(a)).
\end{equation}
\end{enumerate}
\end{theorem}

Although the above definition of Maslov index allows for paths with endpoints not in $\Lambda^0(L_0)$, we shall only be
interested in the case where $\ell(a),\ell(b)\in\Lambda^0(L_0)$ so that $\mu_{L_0}(\ell)$ will in fact be an integer.
The Maslov index $\mu_{L_0}(\ell)$ is a sort of counting of intersections of the path $\ell$ with the {\sl Maslov cycle} (with vertex 
$L_0$) $\Sigma(L_0)=\left\{ L\in \Lambda(V) \; : \; 
L\cap L_0\neq\{0\} \right\}$. If $t_0$ is an isolated instant of intersection, the {\sl Maslov index across} $t_0$, which we denote by
$\mu_{L_0}(\ell,t_0)$, is defined as
\begin{equation}\nonumber
\mu_{L_0}(\ell,t_0) = \mu_{L_0}\left( \ell|_{[t_0-\epsilon,t_0+\epsilon]} \right),
\end{equation}
for $\epsilon$ sufficiently small; the properties of $\mu_{L_0}$ guarantee that $\mu_{L_0}(\ell,t_0)$ is well-defined. 
Let now $\ell:[a,b]\rightarrow\Lambda(V)$ be of class $C^1$. Each velocity vector $\dot{\ell}(t)=d\ell/dt$ identifies 
with a symmetric bilinear form on $\ell(t)$ as follows. Given $t_0$ and $u,v\in\ell(t_0)$, let $u:(t_0-\epsilon, t_0+\epsilon)\rightarrow V$ 
be of class $C^1$
with $u(t)\in \ell(t)$ for all $t$ and $u(t_0)=u$. Then one sets 
\begin{equation}\nonumber
\dot{\ell}(t_0)[u,v]=\omega\left[ \dot{u}(t_0),v\right].
\end{equation}  
An instant $t_0$ of intersection of $\ell$ with $\Sigma(L_0)$ is said to be 
{\sl non-degenerate} if the quadratic form $\dot{\ell}(t_0)$ is non-degenerate on $\ell(t_0)\cap L_0$. We define
the {\sl signature} of such an instant, ${\rm sig}(\ell,L_0,t_0)$, 
as the signature
of the restriction of $\dot{\ell}(t_0)$ to $\ell(t_0)\cap L_0$. In a vicinity of a non-degenerate instant of intersection with $\Sigma(L_0)$
there are no other instants of intersections, and one shows that 
\[
\mu_{L_0}(\ell,t_0)={\rm sig}(\ell, L_0,t_0).
\] 
Thus, for $\ell$ having 
only non-degenerate instants of intersection with $\Sigma(L_0)$, and with endpoints
in $\Lambda^0(L_0)$, the Maslov index is the sum of the corresponding signatures: 
\begin{equation}\nonumber
\mu_{L_0}(\ell) = \sum_{t\in [a,b]} {\rm sig}(\ell, L_0, t).
\end{equation}

\subsection{The Kashiwara index of a triplet of Lagrangians.}\label{sectionkashiwara}

\noindent Given a continuous path $\ell:[a,b]\rightarrow\Lambda(V)$ and two Lagrangians $L_0,L_1\in\Lambda(V)$, it turns out
that the difference of Maslov indices $\mu_{L_1}(\ell) - \mu_{L_0}(\ell)$ only depends on the quadruplet $(L_0,L_1,L_2,L_3)$, where
$L_2=\ell(a)$ and $L_3=\ell(b)$. That number, known in the literature as the {\sl Hörmander four-fold index} of
$(L_0,L_1,L_2,L_3)$,
can be described as half the difference of two {\sl Kashiwara indices}: 
given a triplet $(L_1',L_2',L_3')$ of Lagrangians, 
its Kashiwara index (called in \cite{vergne} the Maslov index of the triplet) $\tau(L_1',L_2',L_3')$ is the
integer defined as the signature of the quadratic form $\srQ$ on $L_1\oplus L_2\oplus L_3$ given by
\begin{equation}\nonumber
\srQ \left[ (v_1,v_2,v_3) \right] = \omega[v_1,v_2] + \omega[v_2,v_3] + \omega[v_3,v_1].
\end{equation}  
Then one has that
\begin{equation}\label{kashiwaraformula}
\mu_{L_1}(\ell) - \mu_{L_0}(\ell) = \frac{1}{2}\tau(L_1,L_0,L_3) - \frac{1}{2}\tau(L_1,L_0,L_2).
\end{equation}
The above equality is implicit in page 181 of \cite{cappellmillerlee} when the path $\ell$ has endpoints transverse to both
$L_0$ and $L_1$, and is stated in \cite{partialsignatures} in the non-transversal case but with a definition of the Maslov index
that differs from ours (we should observe, however, that the claimed relation in \cite{partialsignatures} does not seem to work
since the authors' definition of $\mu_{L_0}(\ell)$ does not change the sign when one changes $\omega$ by $-\omega$).

\medskip

\noindent {\sl Proof of (\ref{kashiwaraformula}) (sketch).}
Very briefly, the proof of (\ref{kashiwaraformula}) is as follows: first one observes that, regarding the right-hand 
side of (\ref{eq2931}), one has ${\rm sig}\,\varphi_{L_0,L_1}(\ell(a))=\tau(L_0,L_1,\ell(a))$ and ${\rm sig}\,\varphi_{L_0,L_1}(\ell(b))=\tau(L_0,L_1,\ell(b))$; this is basically the assertion of Lemma 1.5.4 of \cite{vergne}. Then, given $\ell$ and $L_0$ and $L_1$,
one finds a partition $a=t_0<t_1<\cdots<t_N=b$ of $[a,b]$ and Lagrangians $L_1',...,L_N'\in\Lambda^0(L_0)\cap \Lambda^0(L_1)$ such that
$\ell([t_{i-1},t_i])\subset \Lambda^0(L_i')$ for all $i$. From the properties of the Maslov index stated in
Theorem \ref{theoremmaslovindexdefinition} and the above observation, one obtains
\begin{equation}\nonumber
\mu_{L_1}(\ell) - \mu_{L_0}(\ell) = \frac{1}{2}\sum_{i=1}^N \left( \tau(L_1,L_i',\ell(t_i)) - \tau(L_1,L_i',\ell(t_{i-1}))\right)
- \frac{1}{2}\sum_{i=1}^N \left( \tau(L_0,L_i',\ell(t_{i})) - \tau(L_0,L_i',\ell(t_{i-1}))\right).
\end{equation}
The equality in (\ref{kashiwaraformula}) is then obtained from the above by exploring the chain condition satisfied by $\tau$ (Proposition 1.5.8 of
\cite{vergne}) as well as the antisymmetry of $\tau$.\qed

\section{The spectral flow and bifurcation}\label{sectionspectralflow}

\noindent Throughout this section, $H$ shall denote a separable infinite dimensional real Hilbert space. We shall use 
superscript $\perp$ to denote the orthogonal complement of a subspace and, given a quadratic form $Q$ on $H$, the 
superscript $\perp_Q$ shall denote the orthogonal complement with respect to $Q$.

\subsection{Generalities on the spectral flow.}\label{subsectionspectralflowgeneralities}

\noindent Let $\hat{\gF}=\hat{\gF}(H)$ be the space of self-adjoint Fredholm operators 
on $H$. The space $\hat{\mathfrak{F}}$ has three connected components, $\hat{\gF}_+$, $\hat{\gF}_-$ and $\hat{\gF}_*$. 
The components $\hat{\gF}_{\pm}$ consist of the operators in $\hat{\gF}$ that are essentially positive, resp.\ essentially negative, 
(i.e., that have finite dimensional negative, resp.\ positive, spectral subspaces) and are contractible. As for $\hat{\gF}_*$, it has infinite cyclic 
fundamental group; operators in $\hat{\gF}_*$ are called strongly indefinite.   

Let $L=\{L_s\}_{s\in [a,b]}$ be a continuous path in $\hat{\gF}$; equivalently, via Riez representation, we can consider 
the corresponding path $\cQ=\{\cQ_s\}_{s\in [a,b]}$ of   
quadratic forms of Fredholm type on $H$. The spectral flow of $L$, or of $\cQ$, is the integer ${\rm sf}(L,[a,b])$, or ${\rm sf}(\cQ,[a,b])$,
giving the net number of eigenvalues of $L$ that cross 0 in the 
positive direction as $s$ runs from $a$ to $b$. While for paths in $\hat{\gF}_+$ or in $\hat{\gF}_-$ this definition
can trivially be made precise, for paths in $\hat{\gF}_*$ it requires some effort. There are several approaches in the literature 
that produce the {\sl same} definition of ${\rm sf}(L,[a,b])$ if $L$ has invertible endpoints $L_a$ and $L_b$ (which is the case
we shall be considering); see, for instance, \cite[Sec. 7]{patodi}, 
\cite{phillips}, \cite{recht} and \cite{spectralflowbook}. We shall content ourselves with listing here a few properties of ${\rm sf}(L,[a,b])$ that we shall make use of.

\medskip

\begin{enumerate}[left=4pt .. \parindent]
\item\label{invertibility} ({\sl Normalization}) If $L_s$ is invertible for all $s$, then ${\rm sf}(L,[a,b])=0$;
\item\label{essentiallypositive} if $L$ is a path in $\hat{\gF}_+$ with invertible endpoints, then 
\begin{equation}\nonumber
{\rm sf}(L,[a,b]) = {\rm ind}\, L_a - {\rm ind}\, L_b;
\end{equation}
\item\label{concatenationadditivity} ({\sl Additivity under concatenation}) if $L_1 * L_2$ is the concatenation of two paths $L_1=\{(L_1)_s\}_{s\in [a,b]}$ and $L_2=\{(L_2)_s\}_{s\in[b,c]}$ that
agree for $s=b$, then
\begin{equation}\nonumber
{\rm sf}(L_1*L_2, [a,c]) = {\rm sf}(L_1,[a,b]) + {\rm sf}(L_2,[b,c]);
\end{equation}
\item\label{directsumadditivity} ({\sl Additivity under direct sum}) if $L_1$ and $L_2$ are two paths of self-adjoint Fredholm operators on Hilbert spaces $H_1$ and $H_2$, then 
\begin{equation}\nonumber
{\rm sf}(L_1\oplus L_2,[a,b]) = {\rm sf}(L_1,[a,b]) + {\rm sf}(L_2,[a,b]);
\end{equation}
\item ({\sl Invariance under cogredience}) given a path $L$ in $\hat{\gF}(H)$ and a path $M=\{M_s\}$ in $GL(H)$, then
\begin{equation}\nonumber
{\rm sf}(L,[a,b]) = {\rm sf}(M^*LM,[a,b]);
\end{equation}
\item\label{compactperturbation} ({\sl Compact perturbations of closed paths}) for {\sl closed} paths, ${\rm sf}(L,[a,b])$ is stable under compact perturbations: 
if $L$ is closed and $K=\{K_s\}_{s\in [a,b]}$ is a
closed continuous path of compact self-adjoint operators, then 
\begin{equation}\nonumber
{\rm sf}(L,[a,b]) = {\rm sf}(L+K,[a,b]).
\end{equation}
\end{enumerate}

\subsection{The spectral flow of a restriction to a subspace.}\label{sectionrestrictionspectralflow}

\noindent Let $\cQ=\{\cQ_s\}_{s\in [a,b]}$ be a continuous path of quadratic forms of Fredholm type on $H$. Let $\cV\subset H$ be a closed finite 
codimensional subspace and denote by $\overline{\cQ}=\{\overline{\cQ}_s\}_{s\in [a,b]}$ the path of quadratic forms
on $\cV$ obtained by restricting the $\cQ_s$ to $\cV$. It is clear that $\overline{\cQ}_s$ are also of Fredholm type. 
We shall compute the difference ${\rm sf}(\cQ,[a,b]) - {\rm sf}(\overline{\cQ},[a,b])$ in terms of the endpoints data.
The key observation is the following.
\begin{lemma}\label{lemmaclosedpath}
Suppose that the path $\cQ$ is closed. Then ${\rm sf}(\cQ,[a,b]) = {\rm sf}(\overline{\cQ},[a,b])$.
\end{lemma}
\Proof
Let $L_s$ be the self-adjoint operator representing $\cQ_s$. With respect to the orthogonal decomposition $H=\cV\oplus\cV^\perp$, 
let us write $L_s$ in block form as
\begin{equation}\nonumber
L_s = \begin{pmatrix} A_s & B_s \\ B_s^* & C_s \end{pmatrix} = 
\begin{pmatrix} A_s & {\rm O} \\ {\rm O} & {\rm O} \end{pmatrix} +
\begin{pmatrix} {\rm O} & B_s \\ B_s^* & C_s  \end{pmatrix} = N_s + K_s,
\end{equation}
for $N_s$ and $K_s$ being the first and second summands above, 
respectively. Since $\cV^\perp$ has finite dimension, the operator 
$K_s$ is compact. Thus, by \ref{compactperturbation}. of Sec. \ref{subsectionspectralflowgeneralities}, 
${\rm sf}(L,[a,b]) = {\rm sf}(N,[a,b])$. Also, by \ref{directsumadditivity}. we have ${\rm sf}(N,[a,b])={\rm sf}(A,[a,b])$. 
To conclude, just observe that $A_s$ is the self-adjoint operator representing $\overline{\cQ}_s$.  
\qed

\begin{lemma}\label{lemmaconnected}
The space $\hat{\gF}_*\cap GL(H)$ of invertible strongly indefinite self-adjoint operators is path-connected.
\end{lemma}
\Proof
Fix an orthogonal decomposition $H=H_+\oplus H_-$ with $H_\pm$ infinite dimensional closed subspaces, and let
$\gI$ be the symmetry on $H$ given by $\gI(v_+ + v_-)=v_+ - v_-$.  
Given $S\in\hat{\gF}_*\cap GL(H)$, let us decompose $H$ as the orthogonal sum $H=E_+\oplus E_-$
of the positive and negative spectral subspaces
of $S$. Let $P_\pm$ be the projections onto $E_\pm$ and let $\gI'$ be the symmetry $\gI'=P_+ - P_-$.
We first connect $S$ to $\gI'$ through the path $[0,1] \ni t\mapsto (t I + (1-t)S)\circ P_+ + (-t I + (1-t)S)\circ P_-$.
Then we connect $\gI'$ to $\gI$ as follows. Since $H_\pm$ and $E_\pm$ are infinite dimensional separable Hilbert spaces,
we can find orthogonal transformations $f_\pm : E_\pm \r H_\pm$. As the orthogonal group $O(H)$ is connected, we connect
$I$ to $f_+ \oplus f_-$ through a path $[0,1] \ni t \mapsto A_t \in O(H)$. For each $t$, let $P_\pm^t$ be the projections
onto $A_t E_\pm$. Then $t\mapsto P_+^t - P_-^t$ is a path of symmetries connecting $\gI'$ to $\gI$. 
\qed

\begin{lemma}\label{lemmadirectsum}
Let $s$ be such that $\overline{\cQ}_s$ is non-degenerate. Then $H=\cV \oplus \cV^{\perp_{\cQ_s}}$. 
\end{lemma}
\Proof
We already have $\cV\cap \cV^{\perp_{\cQ_s}}=\{0\}$ since $\overline{\cQ}_s$ is non-degenerate. 
By passing to the quotient $H/{\rm ker}\, \cQ_s$, we can also suppose that $\cQ_s$ is non-degenerate. Thus
the musical map $\cQ_s^\sharp:H \r H^*$ is an isomorphism that sends 
$\cV^{\perp_{\cQ_s}}$ onto the annihilator $\cV^0$ of $\cV$. Therefore 
${\rm dim}\, \cV^{\perp_{\cQ_s}} = {\rm dim}\,\cV^0 = {\rm codim}\,\cV=n$. It follows that 
$(\cV\oplus\cV^{\perp_{\cQ_s}})/\cV \cong \cV^{\perp_{\cQ_s}}$ is a subspace of $H/\cV$ with dimension $n$. 
Hence $\cV\oplus\cV^{\perp_{\cQ_s}}=H$.
\qed

\begin{proposition}\label{propspectralflowdifference}
Let $\cQ$, $\cV$ and $\overline{\cQ}$ be as in the beginning of this section. Suppose that $\cQ$ and $\overline{\cQ}$
have non-degenerate endpoints. Then,
\begin{equation}\label{eqspectralflowrestriction}
{\rm sf}(\cQ,[a,b]) - {\rm sf}(\overline{\cQ},[a,b]) = {\rm ind}\,\cQ_a\big|_{\cV^{\perp_{\cQ_a}}} 
- {\rm ind}\,\cQ_b\big|_{\cV^{\perp_{\cQ_b}}}.
\end{equation}
\end{proposition}
\Proof
We shall only consider the case where the $\cQ_s$ are strongly indefinite, the other cases being simpler.
In this case, it is clear that $\overline{\cQ}_s$ are also strongly indefinite.
By the previous lemma, we have $H=\cV\oplus\cV^{\perp_{\cQ_s}}$ for $s=a,b$. So,
for $s=a,b$, by rotating $\cV^{\perp_{\cQ_s}}$ to $\cV^\perp$ while keeping $\cV$ fixed, we can connect $\cQ_s$, via a path of 
non-degenerate quadratic forms, to the quadratic form $\hat{\cQ}_s$ 
whose block decomposition relative to $H=\cV\oplus \cV^\perp$ is 
\begin{equation}\nonumber
\hat{\cQ}_s =
\begin{pmatrix}
\overline{\cQ}_s & {\rm O} \\
{\rm O} & \cQ_s^\perp
\end{pmatrix}, \quad s=a,b,
\end{equation}
for some quadratic form $\cQ_s^\perp$ on $\cV^\perp$ which is isometric to $\cQ_s|_{\cV^{\perp_{\cQ_s}}}$. 
Next we connect $\hat{\cQ}_b$ to 
\begin{equation}\nonumber
\hat{\cQ}_b'= \begin{pmatrix} \overline{\cQ}_a & {\rm O} \\ {\rm O} & \cQ^\perp_b \end{pmatrix},
\end{equation} 
via a path of non-degenerate quadratic forms on $H$, by connecting $\overline{\cQ}_b$ to $\overline{\cQ}_a$ via a path
of non-degenerate quadratic forms on $\cV$ while keeping the other blocks unaltered; this is possible because of Lemma \ref{lemmaconnected}.  
At last, we connect $\hat{\cQ}_b'$ to 
$\hat{\cQ}_a$ by connecting, in any way, $\cQ_b^\perp$ to $\cQ_a^\perp$ while keeping the other blocks unaltered. 
We end up with a concatenation of paths, $\cQ_a \r \cQ_b \r \hat{\cQ}_b \r \hat{\cQ}_b'\r \hat{\cQ}_a \r \cQ_a$, defined, respectively, on intervals
$[a,b]$, $[b,c]$, $[c,d]$, $[d,e]$, $[e,f]$, resulting in a closed path $\cQ' = \{\cQ_s'\}_{s\in[a,f]}$. Let $\overline{\cQ'}$ be the restriction of $\cQ'$
to $\cV$. By Lemma \ref{lemmaclosedpath}, we have ${\rm sf}(\cQ',[a,f]) = {\rm sf}(\overline{\cQ'},[a,f])$. On the other hand,

\medskip

\noindent 1. since $\cQ'_s$ and $\overline{\cQ'}_s$  are non-degenerate for, respectively, $s\in[b,d]\cup [e,f]$ and $s\in[b,f]$, then
${\rm sf}(\cQ',[b,d])={\rm sf}(\cQ',[e,f])=0$ and ${\rm sf}(\overline{\cQ'},[b,f])=0$;
 
\medskip

\noindent 2. by additivity under orthogonal direct sum (\ref{directsumadditivity}.\ of Sec.\ \ref{subsectionspectralflowgeneralities}), and the fact that 
${\rm dim}\,\cV^\perp <\infty$, we obtain
${\rm sf}(\cQ',[d,e]) = {\rm ind}\,\cQ_b^\perp - {\rm ind}\,\cQ_a^\perp$, which in turn is equal to minus the right-hand side of Eq.\ (\ref{eqspectralflowrestriction})
since  $\cQ_s^\perp$ is isometric to $\cQ_s|_{\cV^{\perp_{\cQ_s}}}$ for $s=a,b$.

\medskip

\noindent The equality (\ref{eqspectralflowrestriction}) follows now from additivity of the spectral flow under concatenation
\qed

\begin{remark}\label{remarkbenevieri}
The computation of the difference ${\rm sf}(\cQ,[a,b]) - {\rm sf}(\overline{\cQ},[a,b])$ was also carried out in the work \cite{benevieri} under the
hypothesis that the self-adjoint representation of $\cQ$  is of the form $L_s=\gI + K_s$ for some (fixed) symmetry $\gI$ and some
path $K$ of self-adjoint compact operators. In their formula, degenerate endpoints were also allowed. 
\end{remark}

\subsection{Abstract variational bifurcation.}\label{sectionabstractbifurcation}

\noindent Let $\cF_s: \cN \subset H \r \R$ be a family of $C^2-$functionals
defined on a neighborhood $\cN$ of the origin, all of which having the origin as a critical point, and
depending smoothly (i.e. $C^2$) on the parameter $s\in [a,b]$. 
A point $s_0\in (a,b)$ is said to be a {\sl bifurcation instant} for the family $\{\cF_s\}$ if there exist  
sequences $(s_n)_n$ in $[a,b]$ and $(u_n)_n$ in $\cN$, with $n\geq 1$, such that

\begin{enumerate}[left=4pt .. \parindent]
\item $u_n$ is a critical point of $\cF_{s_n}$ and $u_n\neq 0$ for all $n$;
\item $s_n \r s_0$ and $u_n \r 0$ as $n\r \infty$.
\end{enumerate}

\noindent Suppose further that, for each $s$, the second derivative at the origin,
\[
\cQ_s=(D^2 \cF_s)_0 : H\times H \r \R,
\]
is a quadratic form of Fredholm type. 
As an application of the Implicit Function Theorem, one shows that if $s_0$ is a bifurcation instant, then $s_0$ is a degenerate instant
of the family $\cQ=\{\cQ_s\}_{s\in[a,b]}$, 
i.e., an instant for which the form $\cQ_{s_0}$ is degenerate; see, for instance, Theorem 8.2.1 of \cite{jost}.

The main result of \cite{recht} is the following sufficient condition for existence of bifurcation (see Theorem 1 of 
\cite{recht}).

\begin{theorem}\label{thmabstractbifurcation}
In the above setting, suppose that the endpoints $\cQ_a$ and $\cQ_b$ are non-degenerate. If ${\rm sf}(\cQ,[a,b])\neq 0$,
then there exists a bifurcation instant in $(a,b)$ for the family $\{\cF_s\}$.
\end{theorem}

As in \cite{recht}, let ${\rm sf}(\cQ,s_0)$ be the spectral flow of $\cQ$ across an isolated degenerate instant of $\cQ$. This is 
defined as the spectral flow ${\rm sf}(\cQ,[s_0-\epsilon, s_0+\epsilon])$ for $\epsilon>0$ such that there is no degenerate instants
in $[s_0-\epsilon, s_0+\epsilon]$ besides $s_0$. Then a corollary of the above theorem is the following. 

\begin{corollary}\label{corollaryabstractbifurcationacross}
An isolated degenerate instant $s_0$ for $\cQ$ will be a bifurcation instant for the family $\{\cF_s\}$ provided that ${\rm sf}(\cQ,s_0)\neq 0$.
\end{corollary}

\section{The variational setup}\label{sectionvariationalsetup}

\noindent In this section we shall establish the (local) variational principle for solutions of 
Eq.\ (\ref{lorentzforceequation}) under the energy constraint $E(\gamma,\dot{\gamma})\equiv \kappa$
and compute the energy-constrained
index form of an electromagnetic geodesic $\gamma$. 

It is important to bear in mind that we shall assume the following dictionary throughout this work: 

\begin{equation*}\label{dictionary}
\begin{aligned}
\mbox{\sl curve} \;\, \gamma:[0,T] \r M & \quad \longleftrightarrow 
\quad (x,T), \;\, \mbox{\sl with} \;\, x:[0,1]\r M, \;\, x(t)=\gamma(tT) \\
\mbox{\sl vector field $\hat{V}$ along $\gamma$} & \quad \longleftrightarrow \quad 
\mbox{\sl vector field $V$ along $x$}, \;\, V(t)=\hat{V}(tT)
\end{aligned}
\end{equation*}

\subsection{Hilbert manifold of paths.}

\noindent Let $p$ and $q$ be two fixed points on $M$ (assumed to be connected).
With the help of an auxiliary complete Riemannian metric on $M$ one can form the space $\Omega_{p,q}([0,1])$ 
of all curves $x:[0,1]\r M$ connecting $p$ to $q$ that have Sobolev regularity $H^1=W^{1,2}$. It is well-known
that $\Omega_{p,q}([0,1])$ has the structure of a Hilbert manifold modeled on the Hilbert space $H^1_0([0,1],\R^n)$, 
and that each tangent space $T_x \Omega_{p,q}([0,1])$ identifies to the space of $H^1-$regular vector fields
$V$ along $x$ vanishing at $t=0$ and $t=1$. 

Following Sec.\ 2.1 of \cite{abbondandolo}, a local parametrization around a {\sl smooth} $x_*\in\Omega_{p,q}([0,1])$ can be constructed
as follows: fix a smooth frame $E_1$, ..., $E_n$ along $x_*$ and, with the help of the auxiliary Riemannian metric on $M$, define
the map 
\begin{equation}\label{mapparametrization}
\psi : [0,1]\times B_\rho^n \r M, \quad \psi(t, \bx) = {\rm exp}_{x_*(t)} \Big( \sum_i {\rm x}_i E_i(t) \Big),
\end{equation}
where $B_\rho^n$ denotes the open ball in $\R^n$ of radius $\rho$ and center ${\bf 0}=(0,...,0)^T$ and $\bx=({\rm x}_1,...,{\rm x}_n)^T$.
We choose $\rho$ small enough so that $\psi(t , \cdot): B^n_\rho \r M$ is a diffeomorphism onto an open set for each $t$. 
By the continuous inclusion $H_0^1 \hookrightarrow C^0$, the subset $H_0^1([0,1],B_\rho^n)$ given by the
$\bx=\bx(t)$ such that $\bx([0,1])\subset B^n_\rho$ is open in $H^1_0([0,1],\R^n)$. 
A parametrization for $\Omega_{p,q}([0,1])$ around $x_*$ is then obtained by setting
\begin{equation}\label{parametrization}
\Psi : H^1_0([0,1],B_{\rho}^n) \r \Omega_{p,q}([0,1]), \quad \Psi(\bx)(t) = \psi(t,\bx(t)).
\end{equation}
It is clear that $\Psi({\bf 0})=x_*$ and that the derivative 
$(D\Psi)_{\bf 0}: H_0^1([0,1],\R^n)\r T_{x_*} \Omega_{p,q}([0,1])$ 
is the isomorphism $\bv = ({\rm v}_1,...,{\rm v}_n)^T \mapsto V=\sum_i {\rm v}_i E_i$.

\subsection{The global 1--form $\eta_\kappa$.}

\noindent Following \cite{gabrielle}, for each given value $\kappa$ of the energy we can characterize the solutions 
$\gamma:[0,T]\r M$ of Eq.\ (\ref{lorentzforceequation}) with energy $\kappa$, and connecting $p$ to $q$,
as the zeros of a certain globally defined $1-$form $\eta_\kappa$ on the manifold $\Omega_{p,q}([0,1])\times \R_+$, 
where $\R_+$ is the set of positive real numbers. 

In general, given a Lagrangian function $L:TM\r \R$, its {\sl $\kappa-$action functional} is the function
$\cA_\kappa^L: \Omega_{p,q}([0,1])\times \R_+ \r \R$ defined by
\begin{equation}\nonumber
\cA_\kappa^L(x,T)=
T \int_0^1 (L(x,\dot{x}/T) + \kappa) dt = \int_0^T (L(\gamma,\dot{\gamma})+\kappa)dt.
\end{equation} 
For us, $L$ will stand for the semi-Riemannian Lagrangian, $L(p,v)=\frac{1}{2}g(p)[v,v]$, and we will drop
the superscript $L$.
The $\cA_\kappa$ is a smooth functional (Lemma 3.1 of \cite{abbondandolo}). Also, by means of the electromagnetic field $\sigma$ we
define a $1-$form $\tau^\sigma$ on $\Omega_{p,q}([0,1])$ via
\begin{equation}\nonumber
\tau^\sigma(x)[V] = \int_0^1 \sigma(x)\left[V , \dot{x} \right]dt = 
\int_0^T g(\gamma)\left[ \hat{V} , Y(\gamma)[\dot{\gamma}] \right]dt, \quad V\in T_x \Omega_{p,q}([0,1]).
\end{equation} 
This $1-$form is smooth, as will be clear from the discussion at the end of 
Sec.\ \ref{subsectionlocalfunctional}. The $1-$form $\eta_\kappa$ on $\Omega_{p,q}([0,1])\times \R_+$ is then defined by
\begin{equation}\nonumber
\eta_\kappa = d\cA_\kappa + \pi_1^* \tau^\sigma 
\end{equation}
where $\pi_1:\Omega_{p,q}([0,1])\times\R_+ \r \Omega_{p,q}([0,1])$ is the projection. By a direct computation,
\begin{equation}\nonumber
\eta_\kappa(x,T)\left[(V,B)\right] = - \int_0^T g(\gamma)\left[ \hat{V} , \frac{D}{dt}\dot{\gamma} - Y(\gamma)[\dot{\gamma}] \right] dt 
+ B \left( \kappa - \frac{1}{T} \int_0^T \frac{1}{2}g(\gamma)\left[\dot{\gamma},\dot{\gamma}\right]dt \right),
\end{equation}
for all $(x,T)\in \Omega_{p,q}([0,1])\times\R_+$ and $(V,B)\in T_{(x,T)} \bigl( \Omega_{p,q}([0,1])\times \R_+  \bigr)= T_x \Omega_{p,q}([0,1]) \times \R$. From there one easily sees:

\begin{proposition}\label{propvariationalprinciple}
A curve $\gamma:[0,T]\r M$ connecting $p$ to $q$ is an electromagnectic geodesic with energy $\kappa$ if, and only if,
the $1-$form $\eta_\kappa$ vanishes at $(x,T)$.
\end{proposition}

\subsection{The local magnetic $\kappa-$action functional.}\label{subsectionlocalfunctional}

\noindent What is important to us is that the $1-$form $\tau^\sigma$, and therefore the $1-$form $\eta_\kappa$, possesses a local primitive. 
Given a smooth $x_*\in\Omega_{p,q}([0,1])$, let $\cO\subset \Omega_{p,q}([0,1])$ be the neighborhood of $x_*$ covered by
the parametrization $\Psi$ in (\ref{parametrization}). Following \cite{assenza}, a primitive for $\tau^\sigma$ defined on $\cO$ can be obtained by setting
\begin{equation}\nonumber
\cC_{\sigma,x_*} :\cO \r \R, \quad  \cC_{\sigma,x_*}(x)=\int_{[0,1]\times[0,1]} (c_{x_*,x})^* \sigma,
\end{equation}
where $c_{x_*,x}=c_{x_*,x}(s,t):[0,1]\times[0,1]\r M$ is any map such that $s\mapsto c_{x_*,x}(s , \cdot)$ is a smooth path connecting
$x_*$ to $x$ in $\cO$; that $\cC_{\sigma,x_*}$ is well-defined, smooth and that $d\cC_{\sigma,x_*}=\tau^\sigma$ can be seen from 
the expression on a local parametrization worked out below. 
The corresponding primitive for $\eta_\kappa$ on the neighborhood $\cO\times\R_+$ of $\{x_*\}\times\R_+$ is then
\begin{equation}\nonumber 
\cA_{\kappa, \sigma, x_*}:\cO \times \R_+ \r \R, \quad \cA_{\kappa,\sigma, x_*}(x,T) = \cA_\kappa(x,T) + \cC_{\sigma,x_*}(x).  
\end{equation}
We shall call $\cA_{\kappa,\sigma,x_*}$ the {\sl local magnetic $\kappa-$action functional} (centered at $x_*$). Then, from 
Proposition \ref{propvariationalprinciple},

\begin{corollary}\label{corollarycriticalpoint}
A curve $\gamma:[0,T]\r M$ such that the corresponding $x$ is in $\cO$ is an electromagnetic geodesic with energy $\kappa$ if,
and only if, $(x,T)$ is a critical point of $\cA_{\kappa,\sigma,x_*}$.
\end{corollary}

In order to establish the smoothness of a construction in Sec.\ \ref{subsectionproofbifurcation}, we shall now work out the expression of 
$\cC_{\sigma,x_*}$ in the parametrization $\Psi$. Let $\psi:[0,1]\times B_\rho^n \r M$ be the map in the definition of $\Psi$.  
Since $d\sigma=0$,
there exists a smooth $1-$form $\theta$ on $[0,1]\times B_\rho^n$ such that $\psi^*\sigma = d\theta$. 
Given $\bx\in H_0^1([0,1],B_\rho^n)$ and a smooth
path $[0,1]\ni s \mapsto \bx_s \in H_0^1([0,1],B_\rho^n)$ connecting the origin ${\bf 0}$ to $\bx$, 
let $F=F(s,t)$ be the map $F(s,t)=(t,\bx_s(t))$. 
Since $F(s,0)$ and $F(s,1)$ are constant and $F(0,t)=(t,{\bf 0})$
we obtain that
\begin{equation}\label{eqlocalexpression}
(\Psi^*\cC_{\sigma,x_*})(\bx) = \int_{[0,1]\times[0,1]}F^*d\theta = \int_0^1 \theta(t,\bx)\left[ (1,\dot{\bx}) \right]dt
- \int_0^1 \theta(t,{\bf 0})[(1,{\bf 0})]dt.
\end{equation}
It is clear from this that $\cC_{\sigma,x_*}$ is well-defined and smooth. Observe also that integrating the $1-$form $\Psi^* \tau^\sigma$
along the path $s\mapsto \bx_s$ one obtains $(\Psi^*\cC_{\sigma,x_*})(\bx)$. This shows that $\tau^\sigma$ is smooth with 
$\tau^\sigma = d\cC_{\sigma,x_*}$ in $\cO$.

\subsection{The energy-constrained index form.}

\noindent Let now $(x,T)$ be a critical point of $\cA_{\kappa,\sigma,x_*}:\cO \times \R_+ \r \R$. Recall that the functional 
$\cA_{\kappa,\sigma,x_*}$ will then possess an intrinsic second derivative at $(x,T)$, which is a symmetric bilinear form 
$(D^2\cA_{\kappa,\sigma,x_*})_{(x,T)}$ on the tangent space 
$T_x \Omega_{p,q}([0,1]) \times \R$ whose quadratic form can be computed as
\begin{equation}\label{eq8021}
(D^2 \cA_{\kappa,\sigma,x_*})_{(x,T)} \bigl[(V,B) , (V,B) \bigr] = \frac{d^2}{ds^2}\Big|_{s=0} \cA_{\kappa,\sigma,x_*}(x_s,T_s),
\end{equation}
for $(-\epsilon,\epsilon)\ni s \mapsto (x_s,T_s) \in \cO\times\R_+$ being any smooth map with $(x_0,T_0)=(x,T)$ and 
$\partial x_s\/ /\partial s|_{s=0}=V$ and $\partial T_s/\partial s|_{s=0}=B$. 

Since any two local primitives for $\eta_\kappa$ must differ by a constant, the 
bilinear form $(D^2 \cA_{\kappa,\sigma,x_*})_{(x,T)}$ will not depend on the choices made, but only on the electromagnetic geodesic. We shall call it the {\sl energy-constrained index form} of $(x,T)$ (or of $\gamma$), and shall denote it by $\cQ_{\kappa,\sigma}(x,T)$, or $\cQ_{\kappa,\sigma}(\gamma)$.

\begin{lemma}
The energy-constrained index form of $(x, T)$ is given by
\begin{equation}\label{eqsecondvariation}
\begin{aligned}
\cQ_{\kappa,\sigma}(x,T)\bigl[(V_1,B_1),(V_2,B_2)\bigr]  \; &  = \;    
\int_0^T g(\gamma) \left[ \frac{D}{dt} \hat{V_1}  , \frac{D}{dt}  \hat{V}_2 \right]dt  
- \int_0^T g(\gamma)\left[\hat{V_1} \, , \, R(\hat{V_2}, \dot{\gamma})\dot{\gamma} -  
(\nabla_{\hat{V_2}} Y)[\dot{\gamma}]  \right.  
\\
& \left. \quad - \; Y\left[\frac{D}{dt} \hat{V}_2 \right] \right] dt  +  \frac{B_2}{T}\int_0^T g(\gamma)\left[ \hat{V_1} , Y[\dot{\gamma}] \right]dt + \frac{B_1}{T}\int_0^T g(\gamma)\left[ \hat{V_2} , Y[\dot{\gamma}] \right]dt \\
& \quad + \;   2\kappa \frac{1}{T} B_1B_2. 
\end{aligned} 
\end{equation}
\end{lemma}
\Proof
Let us compute the right-hand side of Eq.\ (\ref{eq8021}). To shorten notation, let $\bT(s,t)=\dot{x}_s(t)$, $\bS(s,t)=\partial x_s(t)/\partial s$.
Since $d\cC_\sigma = \tau^\sigma$, then
\begin{equation}\nonumber
 \frac{d^2}{ds^2}\Big|_{s=0} \cA_{\kappa,\sigma,x_*}(x_s,T_s)  =   \frac{d^2}{ds^2}\Big|_{s=0} \cA_\kappa(x_s,T_s)  + \frac{d}{ds}\Big|_{s=0}\tau^\sigma(x_s)\left[ \bS\right].
\end{equation}
\noindent $1-$ Recalling the definition of $\cA_\kappa$ and using that $g(x)[\dot{x},\dot{x}]\equiv 2\kappa T^2$, since $g(\gamma)[\dot{\gamma},\dot{\gamma}]\equiv 2\kappa $, the second derivative on the 
right-hand side above splits as
\begin{equation}\nonumber
\frac{d^2}{ds^2}\Big|_{s=0} \cA_\kappa(x_s,T_s) = \frac{1}{T}\frac{d^2}{ds^2}\Big|_{s=0} \int_0^1 \frac{1}{2}g(x_s)\left[ \bT , \bT \right]dt 
-\frac{2B}{T} \frac{d}{ds}\Big|_{s=0} \int_0^1 \frac{1}{2}g(x_s)[\bT , \bT]dt + \frac{2B^2}{T} \kappa
\end{equation}
Now, by quite standard calculations, and using that $D\dot{x}/dt = T Y(x)[\dot{x}]$, since $D\dot{\gamma}/dt=Y(\gamma)[\dot{\gamma}]$,
\begin{eqnarray}
\frac{d}{ds}\Big|_{s=0} \int_0^1 \frac{1}{2}g(x_s)[\bT , \bT]dt  & = & -T \int_0^1 g(x) \bigl[ V , Y(x)[\dot{x}] \bigr]dt, \nonumber \\
\frac{d^2}{ds^2}\Big|_{s=0} \int_0^1 \frac{1}{2}g(x_s)\left[ \bT , \bT \right]dt  & = & -T\int_0^1 g(x)\left[  \frac{D}{ds}\Big|_{s=0}\bS \, , \, Y(x)[\dot{x}] \right ]dt 
+ \int_0^1 g(x) \left[  \frac{D}{dt}V , \frac{D}{dt}V  \right] dt  \nonumber \\
& & - \int_0^1 g(x) \left[ R(V,\dot{x})\dot{x} , V \right] dt  \nonumber
\end{eqnarray}

\noindent $2-$ As for the derivative $d/ds|_{s=0}\tau^\sigma(x_s)\left[ \bS\right]$, one computes
\begin{eqnarray*}
\frac{d}{ds}\Big|_{s=0} \tau^\sigma(x_s)[\bS] & = & \int_0^1 \frac{d}{ds}\Big|_{s=0} g(x_s) \bigl[ \bS , Y(x_s)[\bT] \bigr] dt \\
& = & \int_0^1 g(x) \left[ \frac{D}{ds}\Big|_{s=0}\bS \, , \, Y(x)[\dot{x}] \right]dt + \int_0^1 g(x)\left[ V \, , \, (\nabla_V Y)[\dot{x}] + Y(x)\left[\frac{D}{dt}V\right] \right] dt.
\end{eqnarray*}
From the above computations one obtains that $d^2\cA_{\kappa,\sigma,x_*}/ds^2|_{s=0}$ is equal to the rewrite in terms of $x$ and $V$
of the right-hand side 
of Eq.\ (\ref{eqsecondvariation}) for $(V_1,B_1)=(V_2,B_2)=(V,B)$. The general
expression in Eq.\ (\ref{eqsecondvariation}) is then obtained by a polarization argument once one shows that the second
integral in the right-hand side of Eq.\ (\ref{eqsecondvariation}) is symmetric in $\hat{V}_1$ and $\hat{V}_2$. The proof of the latter
involves an integration by parts and a subsequent use of identity (\ref{eqexteriorderivative}) in a similar way as in the proof of
1.\ of Lemma \ref{lemmaoperators} below.
\qed

\subsection{The operators $\cK(t)$ and $\cD_t$. }

\noindent The sum of the first two terms on the right-hand side of Eq.\ (\ref{eqsecondvariation})  
can be rewritten in a more
compact and symmetrical form, which is formally identical to the second variation of the standard geodesic action functional,
via the introduction of operators
$\cK=\cK(t): T_{\gamma(t)}M \r T_{\gamma(t)}M$ and $\cD_t:\mathfrak{X}(\gamma)\r\mathfrak{X}(\gamma)$,
where $\mathfrak{X}(\gamma)$ is the space of smooth vector fields along $\gamma$, defined by
\begin{eqnarray*}
\cK v & = & R(v,\dot{\gamma})\dot{\gamma} + \frac{1}{2}(\nabla_{\dot{\gamma}} Y)[v] - (\nabla_v Y)[\dot{\gamma}] 
-\frac{1}{4}Y^2[v],\\
\cD_t \hat{V} & = & \frac{D}{d t} \hat{V} - \frac{1}{2}Y[\hat{V}].
\end{eqnarray*} 
These operators can be obtained as the {\sl Jacobi operator} and {\sl dynamical derivative}, 
respectively, of P. Foulon \cite{foulon} as applied to
the Euler-Lagrange system of equations (\ref{lorentzforceequation}).
Observe that $\cD_t$ is a first order linear differential operator which satisfies the Leibniz rule 
$\cD_t (f \hat{V}) = f \cD_t \hat{V} + \dot{f} \hat{V}$, for $f=f(t)$ a $C^1-$function.
Furthermore, we have

\begin{lemma}\label{lemmaoperators}
The operators $\cK$ and $\cD_t$ satisfy the following properties: 
\begin{enumerate}[left=4pt .. \parindent]
\item $\cK$ is $g-$symmetric: $g[\cK(t)u, v] = g[u,\cK(t)v]$, for all $u,v\in T_{\gamma(t)}M$;
\item $\cD_t$ is compatible with the metric: $(d/dt)g[\hat{V},\hat{W}]= g[\cD_t \hat{V},\hat{W}] 
+ g[\hat{V},\cD_t\hat{W}]$, for all $\hat{V},\hat{W}\in\mathfrak{X}(\gamma)$.
\end{enumerate}
\end{lemma}
\Proof
The 2.\ follows at once from the metric compatibility of the connection and the fact that $Y$ is $g-$antisymmetric. Let us prove
1. Since $v\mapsto R(v,\dot{\gamma})\dot{\gamma}$ and $Y^2$ are $g-$symmetric, it suffices to show that 
$v\mapsto S_u[v]=(1/2)(\nabla_u Y)[v] - (\nabla_v Y)[u]$ is $g-$symmetric for all $u$.
On the one hand,
\begin{eqnarray*}
g\bigl[ S_u[v] , w \bigr] - g\bigl[ v, S_u[w]\bigr] & = & \frac{1}{2} g\bigl[ (\nabla_u Y)[v] , w \bigr]  - \frac{1}{2} g \bigl[ v , (\nabla_u Y)[w]  \bigr]  \\
& & - g\bigl[ (\nabla_v Y)[u] , w \bigr]  +  g\bigl[ v , (\nabla_w Y)[u] \bigr] \\
& = & 
 g\bigl[ w , (\nabla_u Y)[v] \bigr] + g\bigl[ u , (\nabla_v Y)[w] \bigr] +
g\bigl[ v , (\nabla_w Y)[u] \bigr] ;
\end{eqnarray*}
for the second equality we used that $\nabla_uY$ is $g-$antisymmetric, which follows from the $g-$antisymmetry of $Y$.
On the other hand, an easy computation shows that the exterior derivative of the $2-$form $\sigma[\cdot \, \, , \, \cdot \,] = g[\, \cdot \, , \, Y[\cdot]]$
is given by
\begin{equation}\label{eqexteriorderivative}
d\sigma[u,v,w] = g\bigl[ u , (\nabla_v Y)[w] \bigr] + g\bigl[ v , (\nabla_w Y)[u] \bigr] + g\bigl[ w , (\nabla_u Y)[v] \bigr].
\end{equation}
Therefore, since $d\sigma= 0$ we conclude the $g-$symmetry of $\cK$.   
\qed

It is straightforward to check that the sum of the first two terms on the right-hand side of Eq.\ (\ref{eqsecondvariation}) can be rewritten as the simpler expression
employing $\cD_t$ and $\cK$,
\begin{equation}\nonumber
\int_0^T g(\gamma)\bigl[ \cD_t \hat{V_1} , \cD_t \hat{V_2} \bigr] dt
- \int_0^T g(\gamma) \bigl[ \hat{V_1} , \cK \hat{V}_2 \bigr] dt.
\end{equation}
The resulting expression for $\cQ_{\kappa,\sigma}(x,T)$ is thus
\begin{equation}\label{secondvariationnew}
\begin{aligned}
\cQ_{\kappa,\sigma}(x,T)\bigl[(V_1,B_1),(V_2,B_2)\bigr] \; & = \;  \int_0^T g(\gamma)\bigl[ \cD_t \hat{V_1} , \cD_t \hat{V_2} \bigr] dt
- \int_0^T g(\gamma) \bigl[ \hat{V_1} , \cK \hat{V}_2 \bigr] dt \\
& \quad \;  + \;    \frac{B_2}{T}\int_0^T g(\gamma)\bigl[ \hat{V_1} , Y[\dot{\gamma}] \bigr]dt + \frac{B_1}{T}\int_0^T g(\gamma)\bigl[ \hat{V_2} , Y[\dot{\gamma}] \bigr]dt
+ 2\kappa \frac{1}{T} B_1B_2.
\end{aligned}
\end{equation}

\section{Energy-constrained Jacobi fields}\label{sectionjacobifields}

\noindent Until the end of this paper, let $\kappa$ be a {\sl non null} value of the energy and let $\gamma:[0,T]\r M$ be an electromagnetic geodesic with energy $\kappa$ connecting $p$ to $q$. In the following we shall
describe the kernel of the corresponding energy-constrained index form $\cQ_{\kappa,\sigma}(x,T)$.

Let $(V,B)$ be in the kernel of $\cQ_{\kappa,\sigma}(x,T)$. In Proposition \ref{propfredholm} it will be shown that $V$ is smooth, so let
us assume this. Thus, integrating by parts the first term in the right-hand side of Eq.\ (\ref{eqsecondvariation}) and rearranging things, we obtain that
\begin{equation}\label{eqkernelindexform}
\begin{multlined}
\hspace{-1.2cm} - \int_0^T g(\gamma) \left[ \frac{D^2}{dt^2} \hat{V_1} + R(\hat{V_1} , \dot{\gamma} ) \dot{\gamma} - (\nabla_{\hat{V_1}}Y)[\dot{\gamma}] -
 Y \left[ \frac{D}{dt}\hat{V_1} \right] - \frac{B_1}{T}Y[\dot{\gamma}] \, , \, \hat{V_2}   \right] dt  \;\; + \;\; \\
 \frac{B_2}{T} \left( \int_0^T g(\gamma)\left[ \hat{V_1} , Y[\dot{\gamma}] \right]dt  + 2 \kappa B_1  \right) \;\; = \;\; 0 \hspace{-1.1cm}
 \end{multlined}
\end{equation}
for all $V_2\in T_x\Omega_{p,q}([0,1])$ and all $B_2\in\R$. It follows from the Fundamental Lemma of the Calculus of Variations that 
\begin{empheq}[left=\empheqlbrace]{align}
&  \frac{D^2}{dt^2} \hat{V}_1 + R(\hat{V}_1 , \dot{\gamma} ) \dot{\gamma} - (\nabla_{\hat{V}_1}Y)[\dot{\gamma}] -
Y \left[ \frac{D}{dt}\hat{V}_1 \right] - \frac{B_1}{T}Y[\dot{\gamma}] = 0, \label{jacobi1} \\
& B_1 = - \frac{1}{2\kappa}  \int_0^T g(\gamma)\left[ \hat{V_1} , Y[\dot{\gamma}] \right]dt . \label{eqB1}
\end{empheq}
Observe that, from Eq.\ (\ref{lorentzforceequation}), and since $\hat{V}_1$ vanishes on the end-points of $[0,T]$, we have that
\begin{equation}\label{eqB2}
\int_0^T  g(\gamma)\bigl[ \hat{V_1} , Y[\dot{\gamma}] \bigr]dt = \int_0^T  g(\gamma) \bigl[ \hat{V_1} , D\dot{\gamma}/dt \bigr] dt =  -\int_0^T g(\gamma) \bigl[ D\hat{V}_1/dt , \dot{\gamma} \bigr]dt.
\end{equation}

\begin{lemma}\label{lemmajacobi}
If a smooth vector field $J$ along $\gamma$ satisfies Eq.\ (\ref{jacobi1}) for some $B_1\in\R$, then $g(\gamma)\bigl[ D J /dt , \dot{\gamma} \bigr]$ is constant as a function of $t$.
\end{lemma}

\Proof
Taking the $g-$inner product of both sides of Eq.\ (\ref{jacobi1}) with $\dot{\gamma}$, 
\begin{eqnarray*}\nonumber 
g(\gamma) \left[ \frac{D^2}{dt^2} J , \dot{\gamma} \right] + g(\gamma) \bigl[ R(J,\dot{\gamma})\dot{\gamma} , \dot{\gamma} \bigr]
- g(\gamma) \bigl[ (\nabla_{J} Y)[\dot{\gamma}] , \dot{\gamma} \bigr]
  - g(\gamma) \left[ Y\left[ \frac{D}{dt}J \right]  , \dot{\gamma} \right] 
 - \frac{B_1}{T} g(\gamma) \bigl[ Y[\dot{\gamma}] , \dot{\gamma} \bigr] &  = & 0.
\end{eqnarray*}
The second term on the left-hand side vanishes. Also, since $Y$ is $g-$antisymetric the same is true of
$\nabla_{J}Y$, and therefore the third and fifth terms also vanish. Recalling that $Y[\dot{\gamma}]=D\dot{\gamma}/dt$, we are left with
\begin{equation}\nonumber
0 \; = \;  g(\gamma) \left[ \frac{D^2}{dt^2} J , \dot{\gamma} \right] - g(\gamma) \left[ Y \left[ \frac{D}{dt} J \right] , \dot{\gamma} \right] 
\; = \; g(\gamma) \left[ \frac{D^2}{dt^2} J , \dot{\gamma} \right] + g(\gamma) \left[ \frac{D}{dt} J , Y[\dot{\gamma}]   \right]  
\; = \; \frac{d}{dt} g(\gamma) \left[  \frac{D}{dt}J , \dot{\gamma}  \right],
\end{equation}
which was what we wanted to show.
\qed

\begin{definition}
An {\sl energy-constrained Jacobi field} along an electromagnetic geodesic $\gamma:[0,T]\r M$
with energy $\kappa\neq0$
is a smooth vector field $J$ along $\gamma$ that satisfies the {\sl energy-constrained Jacobi equation}
\begin{equation}\label{jacobiequation1}
\frac{D^2}{dt^2}J + R(J,\dot{\gamma})\dot{\gamma} - (\nabla_J Y)[\dot{\gamma}] - Y\left[ \frac{D}{dt}J \right]
-\frac{1}{2\kappa} g(\gamma(0)) \left[ \frac{D}{dt}J(0) , \dot{\gamma}(0) \right] Y[\dot{\gamma}] \;\; = \;\; 0. 
\end{equation}
We shall denote by $\cJ_{\gamma,\kappa}$ the space of all energy-constrained Jacobi fields along $\gamma$.
\end{definition}

It follows from Lemma \ref{lemmajacobi} that the quantity $g\bigl[DJ/dt , \dot{\gamma} \bigr]$
is constant for $J\in \cJ_{\gamma,\kappa}$; thus, Eq.\ (\ref{jacobiequation1}) could have been written with 
$g\bigl[DJ/dt , \dot{\gamma} \bigr]|_{t=t_0}$
in place of $g\bigl[DJ/dt , \dot{\gamma} \bigr]|_{t=0}$ for any given $t_0\in[0,T]$. 

In the absence of its last term, Eq.\ (\ref{jacobiequation1}) becomes the well-known linearization of Eq.\ (\ref{lorentzforceequation})
along an electromagnetic geodesic. We shall refer to it as the {\sl ordinary} Jacobi equation along $\gamma$, and accordingly we shall
speak of {\sl ordinary} Jacobi fields. Of course, an ordinary Jacobi field $J$ is also an energy-constrained one precisely when
$g\left[ DJ/dt , \dot{\gamma}\right]\equiv 0$.

\begin{example}\label{examplejacobi}
The vector fields  $J_1^{tan}(t)=\dot{\gamma}(t)$ and $J_2^{tan}(t)=t\dot{\gamma}(t)$ clearly solve the energy-constrained Jacobi equation. 
Reciprocally, any energy-constrained Jacobi field $J$ that is everywhere tangent to $\gamma$ is a linear combination of
$J_1^{tan}$ and $J_2^{tan}$.
\end{example}

From equations  (\ref{jacobi1}), (\ref{eqB1}), (\ref{eqB2}) and Lemma \ref{lemmajacobi}, we can state 

\begin{proposition}\label{propkerneljacobi}
The kernel of $\cQ_{\kappa,\sigma}(x,T)$ is formed by the pairs $(V,B)\in T_x\Omega_{p,q}([0,1])\times \R$ such that $J=\hat{V}$
is an energy-constrained Jacobi field along $\gamma$ (necessarily vanishing for $t=0$ and $t=T$) and
$B=(T/2\kappa) g(\gamma)\left[ DJ/dt , \dot{\gamma}  \right]|_{t=t_0}$, for any $t_0\in [0, T]$.
\end{proposition}    

It follows from this proposition that the bilinear form $\cQ_{\kappa,\sigma}(x,T)$ is degenerate (in the sense of having
non-trivial 
kernel)
if, and only if, there exits a non null $J\in\cJ_{\gamma,\kappa}$ that vanishes at both $t=0$ and $t=T$.

\begin{definition}
An instant $t_0\in (0,T]$ for which there exists a non null $J\in \cJ_{\gamma,\kappa}$ that vanishes at both $t=0$ and $t=t_0$
shall be called an {\sl energy-constrained conjugate instant along $\gamma$}. 
For such an instant $t_0$, let $\mathds{J}_{t_0}$
be the space of all $J\in \cJ_{\gamma,\kappa}$ that vanish at $t=0$ and $t=t_0$, and
let $\mathds{J}'[t_0] =\left\{ DJ/dt|_{t=t_0} \; : \; J\in \mathds{J}_{t_0} \right\}$. We shall say that $t_0$ is {\sl non-degenerate}, if
$\mathds{J}'[t_0]$ is a non-degenerate subspace with respect to the symmetric bilinear form $\mathfrak{g}_{t_0}$ on $T_{\dot{\gamma}(t_0)}M$
defined by
\begin{equation}
\mathfrak{g}_{t_0}\bigl[u,v \bigr] = g(\gamma(t_0)) \bigl[ u,v \bigr] - \frac{1}{2\kappa} g(\gamma(t_0)) \bigl[ \dot{\gamma}(t_0) , u \bigr] g(\gamma(t_0)) \bigl[ \dot{\gamma}(t_0) , v \bigr].
\end{equation}
The {\sl signature} of a non-degenerate conjugate instant $t_0$ is defined as the signature of the restriction of $\mathfrak{g}_{t_0}$ to
$\mathds{J}'[t_0]$.  
\end{definition}

Of course, we can define the form $\mathfrak{g}_t$ for an arbitrary $t\in [0,T]$. Observe that $\mathfrak{g}_t$ coincides with 
$g(\gamma(t))$ in the $g-$orthogonal complement $\bigl\langle \dot{\gamma}(t) \bigr\rangle^\perp$ and that it 
has a one-dimensional kernel generated by $\dot{\gamma}(t)$.

We shall also speak of {\sl ordinary} conjugate instants along $\gamma$; these are defined by changing the term {\sl energy-constrained} by {\sl ordinary}
in the above definition.

\begin{proposition}\label{propriemannlorentz}
If $g$ is Riemannian, or if $g$ is Lorentzian and $\gamma$ is time-like (i.e.\ $\kappa<0$), then, for any fixed-energy conjugate instant $t_0$
the form $\mathfrak{g}_{t_0}$ is positive-definite on $\mathds{J}'[t_0]$. 
\end{proposition}
\Proof
In both cases $g(\gamma(t_0))$ is positive definite on $\bigl\langle \dot{\gamma}(t_0) \bigr\rangle^\perp$, hence, by the above remarks,
it suffices to show that 
$\mathds{J}'[t_0]\cap \bigl\langle \dot{\gamma}(t_0) \bigr\rangle =\{0\}$. Given $J\in \mathds{J}_{t_0}$ with
$DJ/dt|_{t=t_0}=\lambda \dot{\gamma}(t_0)$, the energy-constrained Jacobi field $S(t)=\lambda(-t_0+t)\dot{\gamma}(t)$ satisfies 
$S(t_0)=0=J(t_0)$ and $DS/dt|_{t=t_0}=\lambda\dot{\gamma}(t_0)=DJ/dt|_{t=t_0}$. Hence $S=J$, and therefore 
$\lambda=0$ since $J(0)=0$.  
\qed

\begin{example}\label{exampleconjugate}
Consider the setting in the second part of Example \ref{examplebifurcation1} of a uniform magnetic field on $S^2$.
Let $\gamma:[0,2\pi/\omega]\r S^2$ be a solution to Eq.\ (\ref{eqlorentzsphere}) with $\|\dot{\gamma}\|=1$, hence with angular 
velocity $\omega=\sqrt{1+\beta^2}$. We shall look at the energy-constrained and ordinary conjugate instants along $\gamma$.
 Given a smooth vector field $J$ along $\gamma$, let us decompose it as
$J(t)= x(t)\dot{\gamma}(t) + y(t)(\gamma(t)\times \dot{\gamma}(t))$. A straightforward computation shows that the 
energy-constrained Jacobi equation for $J$ becomes the system of equations
\begin{equation}\label{eqjacobiesfera}
\left\{
\begin{aligned}
 & \ddot{x} + \beta\dot{y}  \; = \; 0 \\
& \ddot{y} - \beta\dot{x} + y   \; = \; -\beta\dot{x}(0) - \beta^2 y(0)
\end{aligned}
\right. .
\end{equation}
Given initial conditions $x(0)$, $\dot{x}(0)$, $y(0)$, $\dot{y}(0)$, the corresponding
solution $(x,y)^T$ to the above equations is
\begin{equation}\nonumber
\begin{pmatrix}
x \\ y
\end{pmatrix} =
\begin{pmatrix}
(\beta/\omega)\cos(\omega t) & -(\beta/\omega) \sin(\omega t) & t & 1 \\
\sin(\omega t) & \cos(\omega t) & 0 & 0 
\end{pmatrix}
\begin{pmatrix}
a_1 \\ a_2 \\ a_3 \\ a_4
\end{pmatrix},
\end{equation}
where $a_1=\dot{y}(0)$, $a_2=y(0)$, $a_3=\dot{x}(0)+\beta y(0)$ and $a_4=x(0)-(\beta/\omega)\dot{y}(0)$. It follows
that the energy-constrained conjugate instants along $\gamma$ are the values of $t\in (0,2\pi/\omega]$ for which 
the following matrix
is singular,
\begin{equation}\nonumber
A=\begin{pmatrix}
(\beta/\omega)\cos(\omega t) & -(\beta/\omega) \sin(\omega t) & t & 1 \\
\sin(\omega t) & \cos(\omega t) & 0 & 0 \\
\beta/\omega & 0 & 0 & 1 \\ 
0 & 1 & 0 & 0
\end{pmatrix}.
\end{equation}
Since ${\rm det}(A)=-t\sin(\omega t)$, $t_1=\pi/\omega$ (half a period) and $t_2=2\pi/\omega$ (a whole period) are the only such instants. As for the
ordinary conjugate instants along $\gamma$: the ordinary Jacobi equation is just the homogenized version of
Eq.\ (\ref{eqjacobiesfera}). After solving it, one arrives at the following matrix
\begin{equation}\nonumber
A'=\begin{pmatrix}
(\beta/\omega)\cos(\omega t) & (\beta/\omega) \sin(\omega t) & t & 1 \\
\sin(\omega t) & -\cos(\omega t) & \beta & 0 \\
\beta/\omega & 0 & 0 & 1 \\ 
0 & -1 & \beta & 0
\end{pmatrix}  
\end{equation}
which must be singular in order for $t$ to be an ordinary conjugate instant. We have 
${\rm det}(A')= t\sin(\omega t) + 2(\omega - 1/\omega)(1-\cos(\omega t))$, which, in terms of $s=\omega t/2$,
can be rewritten as $4\sin(s)( s\cos(s) + \beta^2 \sin(s) )$. The zeros of ${\rm det}(A')$ in $(0,2\pi/\omega]$ are
thus $t_2=2\pi/\omega$ together with $2/\omega$ times the solutions to $\tan(s) = -(1/\beta)s$ in $(0,\pi]$. It is
clear that the latter has only one solution, which lies in $(\pi/2 , \pi)$. Therefore, the ordinary conjugate instants
along $\gamma$ are $t_2$ and some $t_3\in (\pi/\omega , 2\pi/\omega)$. According to Corollary \ref{corollarybifurcation}
and Proposition \ref{propriemannlorentz}, $t_1$ and $t_2$ are the energy-constrained bifurcation instants along $\gamma$, 
whereas $t_2$ and $t_3$ are the ordinary bifurcation instants (see Corollary 5.6 of \cite{portaluri_piccione}).    
Observe that in the limit case $\beta\r 0$ the ratio $t_3/t_1$ approaches one. Indeed, this limit case correspondes to the
geodesic flow of $S^2$, in which case there is no distinction between energy-constrained and ordinary conjugacy. On the
other hand, when $\beta\r \infty$ the ratio $t_2/t_3$ approaches one, mimicking the behavior observed in the flat $\R^2$ with
a constant magnetic field in which case $t_2$ is the only ordinary conjugate instant.
\end{example}

\subsection{Energy-constrained versus ordinary Jacobi fields.}\label{subsectionversus}

\noindent As proved in \cite{giambo} with a different terminology, energy-constrained Jacobi fields along $\gamma$ can 
also be obtained through one-parameter variations of
$(x,T)$ through solutions of Eq.\ (\ref{lorentzforceequation}) with the same energy. More precisely, let $(-\epsilon,\epsilon)\ni s\mapsto (x_s,T_s)$ be a smooth
variation of $(x,T)$ such that $\gamma_s(t) = x_s(t/T_s)$ is a solution of Eq.\ (\ref{lorentzforceequation}) with energy $\kappa$ for every $s$.
If $V(t)=dx_s/ds|_{s=0}$
and $B=dT_s/ds|_{s=0}$,
then $J(t)=\hat{V}(t)$ is an energy-constrained Jacobi field along $\gamma$ with 
$g\left[ DJ/dt , \dot{\gamma} \right] \equiv 2\kappa B/T$. On the other hand, the 
variational vector field of the variation $s\mapsto \gamma_s$ must be an ordinary Jacobi field along $\gamma$. By a
simple calculation, $d\gamma_s/ds|_{s=0} = J - (B/T)t\dot{\gamma}$. This way we obtain
a correspondence between energy-constrained and ordinary Jacobi fields. This motivates the following:

\begin{lemma}\label{lemmaenergyconstrainedordinary}
Given a pair $(J,B)$, where $B\in \R$ and $J$ is a solution to Eq.\ (\ref{jacobi1}) with $B$ 
in place of $B_1$, an ordinary Jacobi field $\check{J}$ along $\gamma$ is obtained by setting 
$\check{J} = J-(B/T)t\dot{\gamma}$.
There holds $g\left[D\check{J}/dt , \dot{\gamma}\right]\equiv 0$ if, and only if, $B=(T/2\kappa)g\left[DJ/dt , \dot{\gamma}\right]$,
in which case $J$ is an energy-constrained Jacobi field. In this case, 
the correspondence
$J \mapsto \check{J}$
restricts to an isomorphism between the space of energy-constrained Jacobi fields $J$ that vanish at $t=0$ and $t=T$, and
the space of ordinary Jacobi fields $\check{J}$ satisfying $g\left[D\check{J}/dt , \dot{\gamma}\right]\equiv 0$
that vanish at $t=0$ and are tangent to $\gamma$ at $t=T$. 
\end{lemma}
\Proof
The first two assertions follow from straightforward computations that we shall omit. As for the third assertion,
just observe that the correspondence $\check{J} \mapsto J= \check{J} - (1/2\kappa T)g [ \check{J}(T),\dot{\gamma} (T)] t\dot{\gamma}$ 
is the inverse of the correspondence $J\mapsto \check{J}$ between the referred spaces.\qed

\section{Symplectic property of energy-constrained Jacobi fields and conjugacy}\label{sectionsymplectic}

\noindent 
There is a natural {\sl presymplectic} structure that is preserved by the energy-constrained Jacobi equation. It is defined, for each $t\in [0,T]$, 
as the following antisymmetric 
bilinear form $\hat{\omega}_t$ on $T_{\gamma(t)}M\oplus T_{\gamma(t)}M$,
\begin{equation*}
\hat{\omega}_t \bigl[ (u_1,v_1) , (u_2,v_2) \bigr]   =  \mathfrak{g}_t \bigl[ v_1,u_2 \bigr] 
- \mathfrak{g}_t \bigl[ u_1 , v_2 \bigr]  + \sigma(\gamma)[u_1,u_2].
\end{equation*}

\begin{lemma}\label{lemmakernelomega}
The form $\hat{\omega}_t$ has a two-dimensional kernel generated by $(\dot{\gamma}, Y[\dot{\gamma}])$ and $(0,\dot{\gamma})$. 
\end{lemma}
\Proof
Let $(u_1,v_1)\in{\rm ker}\,\hat{\omega}_t$. From $0=\hat{\omega}_t[(u_1,v_1),(0,v_2)]=-\mathfrak{g}_t[u_1,v_2]$ for all $v_2$, we obtain 
$u_1\in{\rm ker}\, \mathfrak{g}_t=\langle \dot{\gamma}\rangle$. It follows that $g[v_1,\dot{\gamma}]=g[v_1-Y[u_1] , \dot{\gamma}]$, and 
hence
$0=\hat{\omega}_t[(u_1,v_1),(u_2,0)]=-\mathfrak{g}_t[ u_2 , Y[u_1]-v_1 ]$ for all $u_2$. Therefore, 
$Y[u_1]-v_1\in{\rm ker}\, \mathfrak{g}_t=\langle\dot{\gamma}\rangle$ 
\qed

\begin{lemma}\label{propjacobisymplectic}
The quantity $\hat{\omega}_t\left[ \left( J_1 , DJ_1/dt \right) , \left( J_2 , DJ_2/dt \right) \right]$ remains constant along 
$\gamma$ if $J_1$ and $J_2$ 
are two energy-constrained Jacobi fields.
\end{lemma}
\Proof
Computing the derivative in $t$ of $\hat{\omega}_t\left[ \left( J_1 , DJ_1/dt \right) , \left( J_2 , DJ_2/dt \right) \right]$ 
and showing
that it vanishes is quite standard. Besides Eq.\ (\ref{jacobiequation1}), and the fact that 
$g(\gamma)\bigl[DJ_i/dt , \dot{\gamma} \bigr]$ is constant, in the end one uses   
the relation (\ref{eqexteriorderivative}). We shall omit the details.
\qed

Given vectors $u,w\in T_{p=\gamma(0)} M$, let $J_{u,w}$ be the energy-constrained Jacobi field along $\gamma$ such that 
$J_{u,w}(0)=u$ and $DJ_{u,w}/dt |_{t=0}=w$. This sets up an isomorphism, for each $t\in [0, T]$, 
\begin{equation}\label{isomorphism89}
\left\{
\begin{aligned}
& \Psi_t: T_{p}M \oplus T_{p}M  \rightarrow 
T_{\gamma}M \oplus T_{\gamma}M , \\
& \Psi_t (u,w) = \left(J_{u,w}(t) , \frac{D}{dt}J_{u,w}(t) \right).
\end{aligned}
\right.
\end{equation} 
According to Lemma \ref{propjacobisymplectic}, the above isomorphism preserves the presymplectic form $\hat{\omega}_t$.
It is clear that $t$ is an energy-constrained conjugate instant if, and only if, $(\Psi_t)^{-1}(\{0\}\oplus T_\gamma M)$ has
non null intersection with $\{0\}\oplus T_pM$.
Let $J_1^{tan}$ and $J_2^{tan}$ be as in Example \ref{examplejacobi}. Then $\Psi_t$ descends to an isomorphism 
between quotient spaces,
\begin{equation}
\bar{\Psi}_t: \frac{T_{p}M \oplus T_{p}M}{ \left\langle   (J_1^{tan} , DJ_1^{tan}/dt ) , (J_2^{tan} , DJ_2^{tan}/dt ) \right\rangle|_{t=0}  } \rightarrow 
\frac{T_{\gamma}M \oplus T_{\gamma}M}{  \left\langle (J_1^{tan} , DJ_1^t/dt ) , (J_2^{tan} , DJ_2^{tan}/dt ) \right\rangle   } .
\end{equation}

\noindent On the other hand, Lemma \ref{lemmakernelomega} shows that $(J_1^{tan} , DJ_1^{tan}/dt )$ and $(J_2^{tan} , DJ_2^{tan}/dt )$
span the kernel of $\hat{\omega}_t$. Thus
$\hat{\omega}_t$ descends to a symplectic form on the above quotient spaces which is, therefore, preserved by $\bar{\Psi}_t$. Let us show that
$\bar{\Psi}_t$ also recovers the notion of energy-constrained conjugacy. For this, let $\hat{\pi}_t$ be the quotient map from 
$T_\gamma M \oplus T_\gamma M$ to the above quotient.

\begin{lemma}\label{lemmapsi}
A given $t_0$ is an energy-constrained conjugate instant if, and only if, 
$(\bar{\Psi}_{t_0})^{-1}\hat{\pi}_{t_0} \left( \{0\}\oplus T_\gamma M \right)$ has non null intersection
with $\hat{\pi}_0 \left( \{0\} \oplus T_p M \right)$. 
\end{lemma}
\Proof
It suffices to show that if $J\in \cJ_{\gamma,\kappa}$ is tangent to $\gamma$ at $t=0$ and 
$t=t_0$, then there exists a unique $J'\in \cJ_{\gamma,\kappa}$ that vanishes at $t=0$ and $t=t_0$ and such that 
$J-J'$ is everywhere tangent to $\gamma$. But this is a direct consequence of the fact that 
the energy-constrained Jacobi fields everywhere tangent to $\gamma$ are spanned by $J_1^{tan}$ and $J_2^{tan}$.\qed

We shall now show that up to a symplectic change of variables the map $\bar{\Psi}_t$ is nothing but the map $\bar{\Phi}_t$ in (\ref{mapbarphi}). 
For this, recall that
by means of the semi-Riemannian metric $g$ one can construct a canonical identification
\begin{equation}\label{identificationtangentspace}
T_{(\gamma,\dot{\gamma})} TM \cong T_\gamma M \oplus T_\gamma M
\end{equation}
relative to which the derivative $\Phi_t=(D\zeta_t)_{(p,v)}$ of the Hamiltonian flow $\zeta_t$ at $(p,v)=(\gamma,\dot{\gamma})|_{t=0}$ is described by 
$\Phi_t (u,w) = \left( \check{J}_{u,w} , D\check{J}_{u,w}/dt \right)|_t$, where $\check{J}_{u,w}$ is defined analogously to $J_{u,w}$ by
changing {\sl energy-constrained} to {\sl ordinary}. Moreover,  

\medskip 

\begin{enumerate}[left=4pt .. \parindent]
\item the canonical symplectic form $\omega$ of $TM$ becomes 
\begin{equation}\label{symplecticformstandard}
\omega(\gamma,\dot{\gamma})\bigl[ (u_1,v_1) , (u_2,v_2)  \bigr] = g(\gamma)[v_1,u_2] - g(\gamma)[u_1,v_2] + \sigma(\gamma)[u_1,u_2],
\end{equation}
 
\item the Hamiltonian vector field $S$ is given by 
\begin{equation}\label{hamiltonianvectorfield}
S_{(\gamma,\dot{\gamma})} = \left( \dot{\gamma} , Y[\dot{\gamma}] \right) = \left( J_1^{tan} , DJ_1^{tan}/dt \right),
\end{equation}

\item the tangent space to the energy surface $E^{-1}(\kappa)$ is 
\begin{equation}\label{tangentspaceenergysurface}
T_{(\gamma,\dot{\gamma})} E^{-1}(\kappa) = T_\gamma M \oplus \langle \dot{\gamma} \rangle^\perp.
\end{equation}
Also, the Lagrangian subspace $\bar{L}_{(\gamma,\dot{\gamma})}$ in (\ref{lagrangianpath})  
is the image of $\{0\} \oplus \langle \dot{\gamma} \rangle^\perp$ under the quotient map
$\pi_t:T_{(\gamma,\dot{\gamma})}E^{-1}(\kappa) \r T_{(\gamma,\dot{\gamma})}E^{-1}(\kappa)/\langle S\rangle$. 
\end{enumerate}

\medskip

\noindent From (\ref{hamiltonianvectorfield}) we have that the inclusion 
$\iota_t: T_{(\gamma,\dot{\gamma})}E^{-1}(\kappa) \hookrightarrow T_\gamma M \oplus T_\gamma M$descends to a map
\begin{equation}
\bar{\iota}_t: \frac{T_{(\gamma,\dot{\gamma})}E^{-1}(\kappa) }{\left\langle S_{(\gamma,\dot{\gamma})} \right\rangle } \rightarrow 
\frac{T_{\gamma}M \oplus T_{\gamma}M}{  \left\langle (J_1^{tan} , DJ_1^{tan}/dt ) , (J_2^{tan} , DJ_2^{tan}/dt ) \right\rangle }.
\end{equation}

\begin{lemma}\label{lemmaiota}
The map $\bar{\iota}_t$ is a symplectic isomorphism that relates $\bar{\Phi}_t$ and $\bar{\Psi}_t$, i.e., 
$\bar{\iota}_t \circ \bar{\Phi}_t = \bar{\Psi}_t \circ \bar{\iota}_0$.
\end{lemma}
\Proof
Since ordinary Jacobi fields $\check{J}$ with $D\check{J}/dt \perp \dot{\gamma}$ are energy-constrained Jacobi fields,
we see that   
$\iota_t \circ \Phi_t |_{T_{(p,v)}E^{-1}(\kappa)} = \Psi_t \circ \iota_0$. 
Therefore $\bar{\iota}_t \circ \bar{\Phi}_t = \bar{\Psi}_t \circ \bar{\iota}_0$. Also, it is clear that 
$\iota_t$ relates $\hat{\omega}_t$ with the pull-back of $\omega(\gamma)$ to $T_\gamma M \oplus \langle \dot{\gamma}\rangle^\perp$.
To conclude, observe that the inverse of $\bar{\iota}_t$
is the quotient of the map ${\rm id}\oplus\mathfrak{p}_t:T_\gamma M \oplus T_\gamma M \r T_\gamma M \oplus \langle \dot{\gamma} \rangle^\perp$, 
where $\mathfrak{p}_t:T_\gamma M \r \langle \dot{\gamma}\rangle^\perp$ is the 
$g-$orthogonal projection.\qed

Let $\ell(t)$ be the Lagrangian path defined in (\ref{lagrangianpath}). It follows from the lemma above and 
Lemma \ref{lemmapsi} that $t_0$ is an energy-constrained conjugate instant if, and only if, 
$\ell(t_0) \cap \bar{L}_{(p,v)}\neq\{0\}$.
We shall show that the notions of non-degeneracy and signature for energy-constrained conjugacy match the corresponding ones
for $\ell(t)$.
For each $t$, consider the isomorphism    
\begin{equation}\label{eqisomorphism209}
(\bar{\Phi}_t)^{-1} \circ \pi_t
:\{0\} \oplus \langle \dot{\gamma}(t) \rangle^\perp \r \ell(t).
\end{equation}

\begin{lemma}\label{lemmavelocity}
Via the above isomorphism, the quadratic form $\dot{\ell}(t)$ corresponds to the restriction of 
$g(\gamma(t))$ to $\langle \dot{\gamma}(t)\rangle^\perp$.
\end{lemma}
\Proof
Let us fix a smooth frame $\cE=(E_1,...,E_n)$ along $\gamma$ with $E_n=\dot{\gamma}$ and $E_i\in\langle \dot{\gamma}\rangle^\perp$
for $i<n$. A smooth frame $\cA=\cA(t)$ 
for $\ell(t)$ is then obtained by projecting the $(\Phi_t)^{-1}(0,E_i)$ ($i=1,...,n-1$)
to $T_{(p,v)} E^{-1}(\kappa)/\langle S \rangle$.
Let $A=A(t)$ be the 
$2n\times n$ matrix whose columns are the coordinates of the $(\Phi_t)^{-1}(0,E_i)$ ($i=1,...,n$) in the frame
for $T_\gamma M \oplus T_\gamma M$ given by $\cE \oplus \cE$. Let also
$\boldsymbol{\omega}=\boldsymbol{\omega}(t)$ be the matrix of $\omega(\gamma)$ in the frame $\cE\oplus \cE$ and $G=G(t)$
be the matrix of $g(\gamma)$ in the frame $\cE$. From (\ref{symplecticformstandard}) we obtain that
\begin{equation}\label{symplecticformmatrix}
\boldsymbol{\omega} = \begin{pmatrix} P & -G \\ G & {\rm O} \end{pmatrix} \qquad (\mbox{for some matrix $P=P(t)$}).
\end{equation}  
It is clear that
the matrix of $\dot{\ell}(t)$ in the frame $\cA$ is the matrix $\dot{A}^T \boldsymbol{\omega}A$ with
the last row and last column removed. 
Thus all we have to show is that $\dot{A}^T \boldsymbol{\omega}A=G$.  
Let $\bx=(x_1,...,x_n)^T$ be the coordinates induced by $\cE$. In these coordinates, the ordinary Jacobi equation 
takes the form $\ddot{\bx} + L \bx + N \dot{\bx}=0$ for some square matrices
$L=L(t)$ and $N=N(t)$. Equivalently, we can write  
\begin{equation}\nonumber
\dot{\bz} = H \bz, \quad {\rm with} \; 
H = \begin{pmatrix} {\rm O} & I_n \\ -L & -N \end{pmatrix} \; 
{\rm and} \; \bz=\begin{pmatrix} \bx \\ \dot{\bx} \end{pmatrix}.
\end{equation}
Let $X=X(t)$ be the fundamental matrix solution of $\dot{\bz}=H\bz$. Then $A=X^{-1}\begin{pmatrix} {\rm O} \\ I_n \end{pmatrix}$,
and thus $\dot{A}=-X^{-1}H\begin{pmatrix} {\rm O} \\ I_n \end{pmatrix}$. Also,
since $(\Phi_t)^{-1}$ leaves $\omega(\gamma)$ invariant, we have $(X^{-1})^T \boldsymbol{\omega} X^{-1} = \boldsymbol{\omega}$. Therefore,
from (\ref{symplecticformmatrix}) we obtain
\begin{equation}\nonumber
\dot{A}^T \boldsymbol{\omega}A =  -\begin{pmatrix}{\rm O} & I_n \end{pmatrix} H^T (X^{-1})^T 
\boldsymbol{\omega} X^{-1}\begin{pmatrix} {\rm O} \\ I_n \end{pmatrix} = -
\begin{pmatrix}{\rm O} & I_n \end{pmatrix} H^T \boldsymbol{\omega} \begin{pmatrix} {\rm O} \\ I_n \end{pmatrix}=G.
\end{equation}
This completes the proof. \qed

\begin{proposition}\label{propsignaturemaslov}
An energy-constrained conjugate instant $t_0$ is non-degenerate if, and only if, $t_0$ is a non-degenerate instant
of intersection of $\ell(t)$ with the Maslov cycle $\Sigma(\bar{L}_{(p,v)})$. Furthermore, for such an instant we have 
${\rm sig}\left(\ell, \bar{L}_{(p,v)}, t_0 \right) = {\rm sig}\left(\mathfrak{g}_{t_0}|_{\mathds{J}'[t_0]}\right)$, and hence
$\mu_\kappa(\gamma,t_0)={\rm sig}\left(\mathfrak{g}_{t_0}|_{\mathds{J}'[t_0]}\right)$.
\end{proposition}
\Proof
Firstly, it is easily seen that $\mathds{J}'[t_0]\cap \bigl\langle \dot{\gamma}(t_0) \bigr\rangle =\{0\}$, hence the 
$g-$orthogonal projection $\mathfrak{p}_{t_0}:T_{\gamma(t_0)}M\r \bigl\langle \dot{\gamma}(t_0) \bigr\rangle^\perp$ restricts to an injection
$\mathds{J}'[t_0]\r \bigl\langle \dot{\gamma}(t_0) \bigr\rangle^\perp$. Since $\mathfrak{g}_{t_0}$ has kernel 
$\bigl\langle \dot{\gamma}(t_0) \bigr\rangle$ and coincides with $g(t_0)$ on $\bigl\langle \dot{\gamma}(t_0) \bigr\rangle^\perp$, the
injection $\mathfrak{p}_{t_0}|_{\mathds{J}'[t_0]}$ will be isometric if we endow $\mathds{J}'[t_0]$ and $\bigl\langle \dot{\gamma}(t_0) \bigr\rangle^\perp$  with the restrictions of $\mathfrak{g}_{t_0}$ and $g(t_0)$, respectively. Hence, if we show that $\ell(t_0)\cap \bar{L}_{(p,v)}$ corresponds to 
$\mathfrak{p}_{t_0}(\mathds{J}'[t_0])$ under the isomorphism (\ref{eqisomorphism209}),
the result will follow from Lemma \ref{lemmavelocity}. For this, let $\hat{\pi}_t$ and $\bar{\iota}_t$ be as 
in Lemmas \ref{lemmapsi} and \ref{lemmaiota}, respectively. Since $\bar{\iota}_t$ sends  
$\bar{L}_{(\gamma,\dot{\gamma})}$ onto $\hat{\pi}_t\left(\{0\}\oplus T_{\gamma}M\right)$ for each $t$,
it follows from Lemma \ref{lemmaiota} that $\bar{\iota}_{0}$ sends $\bar{L}_{(p,v)}\cap\ell(t_0)$ onto 
$\hat{\pi}_0\left( \{0\} \oplus T_pM \right) \cap (\bar{\Psi}_{t_0})^{-1} \hat{\pi}_{t_0}\left( \{0\}\oplus T_\gamma M \right)$.
On the other hand, from the description of $\Psi_t$ in terms of energy-constrained Jacobi fields, and from Lemma \ref{lemmapsi},
we see that the latter intersection is equal to $(\bar{\Psi}_{t_0})^{-1} \hat{\pi}_{t_0}\left( \{0\}\oplus \mathds{J}'[t_0] \right)$.
It follows that $\bar{\iota}_{t_0}$ sends $\bar{\Phi}_{t_0}\left( \bar{L}_{(p,v)} \cap \ell(t_0) \right)$ onto 
$\hat{\pi}_{t_0}\left( \{0\}\oplus \mathds{J}'[t_0] \right)$. Therefore, since the inverse of $\bar{\iota}_{t_0}$
is the quotient of the map ${\rm id}\oplus\mathfrak{p}_{t_0}:T_\gamma M \oplus T_\gamma M \r T_\gamma M\oplus\bigl\langle \dot{\gamma}(t_0)\bigr\rangle^\perp$, we conclude that $\ell(t_0)\cap \bar{L}_{(p,v)}$ corresponds to 
$\mathfrak{p}_{t_0}(\mathds{J}'[t_0])$ under (\ref{eqisomorphism209}).\qed

Another consequence of Lemma \ref{lemmavelocity} is the following.

\begin{proposition}\label{propisolatedinstant}
There exists $\epsilon>0$ such that there is no energy-constrained conjugate instant in $(0,\epsilon]$. 
\end{proposition} 
\Proof
We have that $\ell(0)\cap \bar{L}_{(p,v)}=\bar{L}_{(p,v)}\cong \langle v \rangle^\perp$ and $\dot{\ell}(0)$ is the restriction of $g(p)$ to 
$\langle v \rangle^\perp$. Since $\kappa\neq 0$, that restriction is non-degenerate. Thus $\ell(t)$ has a non-degenerate intersection with
$\Sigma(\bar{L}_{(p,v)})$ at $t=0$. Since such intersections are isolated, the result follows.
\qed

\section{The energy-constrained spectral flow of an electromagnetic geodesic }\label{sectionanalytic}

\noindent In what follows we shall denote by $\mathds{H}$ the space $H^1_0([0,1],\R^n)\times \R$. On $\mathds{H}$ we choose the Hilbert inner product
\begin{eqnarray*}
\bigl\langle (\bv_1,B_1),(\bv_2,B_2) \bigr\rangle_{\mathds{H}}  & = &  \langle \bv_1 , \bv_2 \rangle_{H^1} + B_1B_2  \\
& = &  
\langle \bv_1 , \bv_2 \rangle_{L^2} + \langle \dot{\bv}_1 , \dot{\bv}_2 \rangle_{L^2} + B_1B_2.
\end{eqnarray*}

By letting the endpoint $q$ vary along $\gamma$ and computing the energy-constrained index form of the corresponding
electromagnetic geodesic, we obtain a one-parameter family of quadratic forms
$\cQ_{\kappa,\sigma}(x_s,s)$, $s\in(0,T]$, with varying domains, 
where $x_s:[0,1]\r M$ is the curve $x_s(t)=\gamma(st)$. We then bring this to a path $\{Q_s\}$ of quadratic forms on $\mathds{H}$ as follows: 
fixed a smooth frame $E_1$, ..., $E_n$ along $\gamma$, the correspondence
\begin{equation}\nonumber
\bv=\bv(t)=({\rm v}_1(t),...,{\rm v}_n(t))^T \quad \mapsto \quad V^s=V^s(t) = \sum_i {\rm v}_i(t) E_i(st)
\end{equation}
is an isomorphism between $H_0^1([0,1],\R^n)$ and the tangent space $T_{x_s}\Omega_{p,\gamma(s)}([0,1])$. We then obtain an isomorphism, for each $s\in (0,T]$,
\begin{equation}\label{isomorphismtrivialization}
\left\{ 
\begin{aligned}
& \mathds{H} \stackrel{\simeq}{\longrightarrow} T_{x_s} \Omega_{p,\gamma(s)}([0,1]) \times \R , \\
& (\bv,B)\mapsto (V^s,sB).
\end{aligned}
\right.
\end{equation}
The quadratic forms $\{Q_s\}$ are then defined as the ones corresponding to the $\cQ_{\kappa,\sigma}(x_s,s)$ via the above isomorphism,
\begin{equation}\label{eqquadraticform}
Q_s \bigl[ (\bv_1,B_1),(\bv_2,B_2) \bigr] = \cQ_{\kappa,\sigma}(x_s,s)\bigl[ (V_1^s , sB_1) , (V_2^s , sB_2) \bigr].
\end{equation}
Any two choices of frames will amount to paths 
$\{Q_s\}$ and $\{Q_s'\}$ that are isometric, i.e., the $Q_s'$ are of the form $M_s^*Q_s$ for some continuous (in fact smooth) path $\{M_s\}$
in $GL(H)$. Thus $\{Q_s\}$ and $\{Q_s'\}$ will have essentially the same properties. 
In order to obtain a simple expression for the $Q_s$
we shall make the assumptions that the frame 
$E_1$, ..., $E_n$ is both $\cD_t-$parallel and $g-$orthonormal, i.e.,
\begin{equation}\label{matrixInp}
\cD_t E_i = 0 \quad {\rm and} \quad \bigl(g(\gamma)[E_i,E_j]\bigr) = I_{n,p}, \quad {\rm for} \quad (I_{n,p})_{ij} = 
\left\{ 
\begin{aligned}
 -1, & \quad {\rm if} \quad  i=j\leq p, \\ 
 1,  & \quad {\rm if} \quad i=j>p, \\ 
 0,  & \quad {\rm if} \quad i\neq j, 
\end{aligned} \right. 
\end{equation}
where $0\leq p \leq n$ is the index of the metric $g$; the existence of such a frame follows from $\cD_t$ being a first order 
linear differential operator which satisfies 2.\ of
Lemma \ref{lemmaoperators}.  
With respect to this frame,
let
\begin{align*}
&\bK =   \bK(t) \; \mbox{be the product of $I_{n,p}$ with the matrix of $\cK(t)$;}\\
&\bxi = \bxi(t) \; \mbox{be the product 
of $I_{n,p}$ with the column-vector} \\
& \qquad \qquad \mbox{of the coordinates of $Y[\dot{\gamma}(t)]$;} \\
& \bK_s(t)  = s^2 \bK(st)  \quad \mbox{and} \quad \bxi_s(t)  =  s^2 \bxi(st).
\end{align*}
Observe that since $\cK$ is $g-$symmetric, then $\bK^T=\bK$. Let also $\langle,\rangle$ denote the standard scalar product of $\R^n$.   

\begin{lemma}
The quadratic forms $Q_s$ are given by
\begin{equation}\label{secondvariationparametric}
\begin{aligned}
 Q_s \bigl[( \bv_1,B_1) , (\bv_2 , B_2) \bigr] \; & = \;  \frac{1}{s}\int_0^1 \langle I_{n,p} \dot{\bv}_1 , \dot{\bv}_2 \rangle dt 
- \frac{1}{s} \int_0^1 \langle \bK_s \bv_1 ,  \bv_2 \rangle dt  \\
& \quad \;  + \; B_1 \frac{1}{s} \int_0^1 \langle \bv_2 , \bxi_s  \rangle dt 
 + B_2 \frac{1}{s} \int_0^1 \langle \bv_1 , \bxi_s \rangle dt +2\kappa s B_1B_2 .
\end{aligned}
\end{equation}
In particular, $\{Q_s\}$ is a smooth path.
\end{lemma}
\Proof
The right-hand side of (\ref{eqquadraticform}) was already computed in Eq.\ (\ref{secondvariationnew}) once one 
replaces $T$ by $s$ and $B_i$ by $sB_i$ in the latter. To show that it matches (\ref{secondvariationparametric}),
observe first that if $V(t)=\sum_i v_i(t)E_i(st)$, then 
$\hat{V}(t)=\sum_i v_i(t/s) E_i(t)$, and then use that $E_1$, ..., $E_n$ is an orthonormal $\cD_t-$parallel frame 
and that $\cD_t$ satisfies the Leibniz rule; we leave the details to the reader.
\qed

\begin{proposition}\label{propfredholm}
The forms $Q_s$ are of Fredholm type. Also, if $(\bv,B)$ is in the kernel of $Q_s$ then $\bv$ is smooth. 
\end{proposition}
\Proof
Let $-d^2/dt^2 : H^2([0,1],\R^n)\cap H_0^1([0,1],\R^n) \subset L^2([0,1],\R^n)\r L^2([0,1],\R^n)$ be the Dirichlet Laplacian.
It is well-known that $1-d^2/dt^2$ has a bounded inverse $(1-d^2/dt^2)^{-1}:L^2 \r H^2\cap H_0^1$ and it is clear that
\[
\bigl\langle \bv_1 , \bv_2 \bigr\rangle_{L^2} = \Bigl\langle (1-d^2/dt^2)^{-1} \bv_1 \, , \, \bv_2 \Bigr\rangle_{H^1}, \quad
{\rm if} \quad \bv_1\in L^2,\; \bv_2\in H^1_0.
\]
Thus we can write, for $\bv_1$ and $\bv_2$ in $H^1_0$,
\begin{multline*}
\frac{1}{s} \int_0^1 \langle I_{n,p} \dot{\bv}_1 , \dot{\bv_2} \rangle dt - \frac{1}{s}\int_0^1 \langle \bK_s \bv_1 , \bv_2 \rangle dt \;\; = \;\;
\frac{1}{s} \bigl\langle I_{n,p} \bv_1 , \bv_2 \bigr\rangle_{H^1}
-\frac{1}{s} \bigl\langle I_{n,p} \bv_1 , \bv_2 \bigr\rangle_{L^2} - \frac{1}{s} \bigl\langle \bK_s \bv_1 , \bv_2 \bigr\rangle_{L^2} \\
 \;\; = \;\; \frac{1}{s} \bigl\langle I_{n,p} \bv_1 , \bv_2 \bigr\rangle_{H^1} - \frac{1}{s} \Bigl\langle (1-d^2/dt^2)^{-1} I_{n,p} \bv_1 , \bv_2 \Bigr\rangle_{H^1}
- \frac{1}{s}\Bigl\langle (1-d^2/dt^2)^{-1}\bK_s \bv_1 , \bv_2 \Bigr\rangle_{H^1}
\end{multline*}
and 
\begin{multline*}
B_1\frac{1}{s} \int_0^1 \langle \bv_2 , \bxi_s  \rangle dt 
+ B_2\frac{1}{s} \int_0^1 \langle \bv_1 , \bxi_s \rangle dt +2\kappa s B_1B_2  \;\; = \;\; 
\frac{1}{s} \bigl\langle B_1\bxi_s , \bv_2 \bigr\rangle_{L^2} + \left( \frac{1}{s} \bigl\langle \bv_1 , \bxi_s \bigr\rangle_{L^2} + 2\kappa s B_1  \right)B_2 \\
 \;\; = \; \; \frac{1}{s} \Bigl\langle B_1(1-d^2/dt^2)^{-1} \bxi_s , \bv_2 \Bigr\rangle_{H^1} + \left( \frac{1}{s} \bigl\langle \bv_1 , \bxi_s \bigr\rangle_{L^2} + 2\kappa s B_1  \right)B_2.
\end{multline*}
Therefore
\begin{equation}\label{eq4923}
\begin{aligned}
Q_s \bigl[ ( \bv_1,B_1) , (\bv_2 , B_2) \bigr] & \;\;=\;\; 
\frac{1}{s} \Bigl\langle I_{n,p} \bv_1 - (1-d^2/dt^2)^{-1}\left( I_{n,p} \bv_1 + \bK_s\bv_1 - B_1\bxi_s \right) \, , \, \bv_2 \Bigr\rangle_{H^1}
\\ &  \quad\quad + \left( \frac{1}{s} \bigl\langle \bv_1 , \bxi_s \bigr\rangle_{L^2} + 2\kappa s B_1  \right)B_2 ,
\end{aligned}
\end{equation}
from where it follows that
the self-adjoint operator $L_s$ which represents $Q_s$ is given by $L_s=\cI-K$, with
\begin{eqnarray*}
\cI(\bv,B) & = & \Bigl( \, \frac{1}{s} I_{n,p} \bv \, , \,  2\kappa s B \, \Bigr), \\ 
K(\bv,B) & = & \frac{1}{s} \Bigl( \, (1-d^2/dt^2)^{-1}( I_{n,p}  + \bK_s )\bv \, , \, 0 \Bigr) 
-\frac{1}{s} \Bigl( B \overline{\bxi}_s \, , \, \bigl\langle \bv , \bxi_s \bigr\rangle_{L^2} \Bigr)
\end{eqnarray*}
and $\overline{\bxi}_s = (1-d^2/dt^2)^{-1} \bxi_s$.
It is clear that $\cI$ is an isomorphism of $\mathds{H}$. As for $K$, it is a compact operator, since the operator 
$(1-d^2/dt^2)^{-1}( I_{n,p} + \bK_s):H_0^1 \r H^1_0$ factors through the compact injection $H^1_0 \hookrightarrow L^2$, and
the operator $\mathds{H}\r \mathds{H}$, 
$(\bv , B)\mapsto \bigl( B \overline{\bxi}_s \, , \, \langle \bv , \bxi_s \rangle_{L^2} \bigr)$ has rank two (unless $Y[\dot{\gamma}(t)]=0$ for $t\leq s$, in which
case it has rank one). Hence $L_s$ is
Fredholm. 

If $(\bv , B)$ is in the kernel of $L_s$, then 
$\bv = I_{n,p} (1-d^2/dt^2)^{-1}( I_{n,p} + \bK_s)\bv - B I_{n,p} \overline{\bxi}_s$.  
This, together with the fact that $(1-d^2/dt^2)^{-1}$ maps $H^l$ into $H^{l+2}$ for each $l=1,2,...$, shows that $\bv$ is smooth. 
\qed

The degenerate instants for $Q=\{Q_s\}$ are the energy-constrained conjugate instants along $\gamma$ and, according to 
Proposition \ref{propisolatedinstant}, there is no such instant in $(0,\epsilon]$ for $\epsilon$ sufficiently small. 
Therefore we can define 
\begin{definition}
Suppose that $T$ is not an energy-constrained conjugate instant along $\gamma$. The {\sl energy-constrained spectral flow}
of $\gamma$, ${\rm sf}_\kappa(\gamma)$, is defined as ${\rm sf}_\kappa (\gamma) = {\rm sf}(Q,[\epsilon, T])$ for $\epsilon>0$ sufficiently small.
\end{definition}

As already observed, a change of frame used to define $Q=\{Q_s\}$ will produce a path $Q'=\{Q_s'\}$ that is isometric to the former; in terms
of the self-adjoint representations $L=\{L_s\}$ and $L'=\{L_s'\}$ of $Q$ and $Q'$, this means that $L$ and $L'$ are cogredient.   
Since the spectral flow is invariant under cogredience, it follows that the above definition does not depend
on the choices made.

\section{Proofs of the main results}\label{sectionproofs}

\noindent Throughout we shall suppose that $T$ is neither an energy-constrained nor an ordinary
conjugate instant along $\gamma$, as well as $\epsilon>0$ is such that there are no such instants in $(0,\epsilon]$. 
We shall continue to assume the setting of Sec.\ \ref{sectionanalytic}.

\subsection{Computing the difference of spectral flows.}

\noindent Let $\cV\subset \mathds{H}$ be the subspace $\cV=H^1_0([0,1],\R^n)\times \{0\}$, and let $\overline{Q}=\{\overline{Q}_s\}$ be the restrictions of $Q_s$ to $\cV$;
the form $\overline{Q}_s$ corresponds to the usual index form of the electromagnetic geodesic $\gamma|_{[0,s]}$. As is well-known, degenerate instants for 
$\{\overline{Q}_s\}$ are the ordinary conjugate instants along $\gamma$. The
spectral flow ${\rm sf}(\overline{Q},[\epsilon, T])$ shall be called the {\sl ordinary spectral flow} of $\gamma$ and denoted by ${\rm sf}(\gamma)$. 

The first step in the proof of Theorem \ref{morseindextheorem} is to compute the difference ${\rm sf}_\kappa(\gamma)-{\rm sf}(\gamma)$. According to 
Proposition \ref{propspectralflowdifference}, we have
\begin{equation}\label{eqdifferencespectralflow}
{\rm sf}_\kappa(\gamma)-{\rm sf}(\gamma) = {\rm ind}\, Q_\epsilon\big|_{\cV^{\perp_{Q_\epsilon}}} -
{\rm ind}\, Q_T \big|_{\cV^{\perp_{Q_T}}}.
\end{equation} 
So let us compute the above indices. First, observe that

\begin{lemma}
Given $s\in(0,T]$, if $(\bv,B)\in \cV^{\perp_{Q_s}}$ then $\bv$ is smooth.
\end{lemma}
\Proof
The proof is identical to the proof we gave for the case $(\bv, B)\in{\rm ker}\, Q_s$ in Proposition \ref{propfredholm}.\qed

\begin{proposition}\label{propindexperp}
Let $s\in (0,T]$ not be an ordinary conjugate instant along $\gamma$. There is a unique ordinary Jacobi field $\check{J}^*_s$ along $\gamma$ 
such that $\check{J}^*_s(0)=0$ and
$\check{J}^*_s(s)=\dot{\gamma}(s)$. In terms of $\check{J}^*_s$, we have
\begin{equation}\nonumber
{\rm ind}\, Q_s\big|_{\cV^{\perp_{Q_s}}} = 
\left\{
\begin{aligned}
& 1,   \quad {\rm if} \quad g(\gamma)\left[ D \check{J}^*_s/dt , \dot{\gamma} \right]\big|_{t=s} <0, \\
& 0,  \quad {\rm if} \quad g(\gamma)\left[ D \check{J}^*_s/dt , \dot{\gamma} \right]\big|_{t=s} \geq 0.
\end{aligned}
\right.
\end{equation} 
Moreover, such an $s$ is not an energy-constrained conjugate instant if, and only if, 
$g(\gamma) \left[ D \check{J}^*_s/dt , \dot{\gamma} \right]\big|_{t=s}\neq 0$.  
\end{proposition}
\Proof
The existence and uniqueness of such $\check{J}^*_s$ follows at once from the assumption 
that $s$ is not an ordinary conjugate instant along $\gamma$. 
Also, since $\cV$ is one-codimensional and $\mathds{H}=\cV\oplus\cV^{\perp_{Q_s}}$ (Lemma \ref{lemmadirectsum}), it follows that 
${\rm dim}\,\cV^{\perp_{Q_s}}=1$. In order to describe $Q_s\big|_{\cV^{\perp_{Q_s}}}$, we shall work directly with $\cQ_{\kappa,\sigma}(x,s)$ instead of 
$Q_s$. So let us take a non null $(V,B)$ in the $\cQ_{\kappa,\sigma}(x,s)-$orthogonal complement of 
$ T_x\Omega_{p,q}([0,1])\times \{0\}$. From the above lemma, $V$ is smooth.
By making $B_2=0$ in Eq.\ (\ref{eqkernelindexform}), we see that $J=\hat{V}$ is a solution to Eq. (\ref{jacobi1}) with $B$ in place of $B_1$
and $s$ in place of $T$.
By making $(V_1,B_1)=(V_2,B_2)=(V,B)$ in the left-hand side of Eq.\ (\ref{eqkernelindexform}), we thus obtain
\begin{equation}\nonumber
\cQ_{\kappa,\sigma}(x,s)[(V,B),(V,B)]  =   \frac{B}{s} \left( \int_0^s g(\gamma)\bigl[ J , Y[\dot{\gamma}] \bigr]dt  + 2 \kappa B \right) .
\end{equation}
Integrating by parts as in (\ref{eqB2}) and using Lemma \ref{lemmajacobi}, we arrive at
\begin{equation}\label{eq27810}
\cQ_{\kappa,\sigma}(x,s)[(V,B),(V,B)] = -B g(\gamma)\left[ D J/dt, \dot{\gamma} \right] \big|_{t=s} + (2\kappa/s)B^2.
\end{equation} 
Now, according to Lemma \ref{lemmaenergyconstrainedordinary}, $\check{J}=J-(B/s)t\dot{\gamma}$ is an ordinary Jacobi field. Rewriting the 
right-hand side of (\ref{eq27810}) in terms of $\check{J}$,  
we obtain
\begin{equation}\label{eq3520}
\cQ_{\kappa,\sigma}(x,s)[(V,B),(V,B)] = \frac{1}{2\kappa}g(\gamma)\left[ \check{J} , \dot{\gamma}\right]
g(\gamma)\bigl[ D\check{J}/dt , \dot{\gamma} \bigr]\, \big|_{t=s},
\end{equation}
since $\check{J}(T)=-B\dot{\gamma}(s)$ (recall that $J(s)=0$). Since the pair $(V,B)$ is non null, and $s$ is not an ordinary conjugate instant,
we must have $B\neq 0$, otherwise $J$ would be a non null ordinary Jacobi field vanishing at $t=0$ and $t=s$. Thus, by dividing by $-B$, we 
could have taken $(V,B)$ with $B=-1$, in which case $\check{J}=\check{J}^*_s$ and the right-hand side of (\ref{eq3520}) would be 
$g(\gamma)\bigl[ D\check{J}^*_s/dt , \dot{\gamma} \bigr]\big|_{t=s}$. This concludes the proof the first assertion. The second assertion
follows from $T_x\Omega_{p,q}([0,1])\times \mathds{R}$ being the direct sum of $T_x\Omega_{p,q}([0,1])\times\{0\}$ with its
$\cQ_{\kappa,\sigma}(x,s)-$orthogonal complement, and from  
$(V,B)$ being a generator of that one-dimensional subspace.\qed

\begin{proposition}\label{propepsilonindex}
It holds that 
\begin{equation}\nonumber
{\rm ind}\,Q_\epsilon\big|_{\cV^{\perp_{Q_\epsilon}}} = 
\left\{
\begin{aligned}
& 0, \quad {\rm if} \quad \kappa>0,\\
& 1, \quad {\rm if} \quad \kappa<0.
\end{aligned}
\right.
\end{equation}
\end{proposition}
\Proof
For $s\in(0,\epsilon]$ let $\hat{Q}_s=sQ_s$. Since $\bK_s(t)=s^2\bK(st)$ and 
$\bxi_s(t)=s^2\bxi(st)$, the quadratic form $\hat{Q}_s$ writes
\begin{equation}\label{hatQ}
\frac{1}{s^2}\hat{Q}_s[(\bv,B)] = \frac{1}{s^2}\alpha[(\bv,B)] + \beta_s[(\bv,B)] + 2\kappa B^2,
\end{equation}
where $\alpha$ and $\beta_s$ are the quadratic forms on $\mathds{H}$ given by
\[
\alpha[(\bv,B)] = \int_0^1\langle I_{n,p}\dot{\bv},\dot{\bv} \rangle dt, \quad
\beta_s [(\bv , B)] = -\int_0^1 \langle \bK(st)\bv(t) , \bv(t) \rangle dt + 2B\int_0^1\langle \bv(t) , \bxi(st) \rangle dt.
\]
It is clear that $\{\beta_s\}_s$ remains bounded for $s\in(0,\epsilon]$, hence $\hat{Q}_s  = \alpha + O(s^2)$ for $s\r 0$. From this
one can conclude that, as $s\r 0$, the one-dimensional subspaces $\cV^{\perp_{\hat{Q}_s}}=\cV^{\perp_{Q_s}}$ converge to 
$\cV^{\perp_\alpha} = \{{\bf 0}\}\times \mathds{R}$ as fast as $s^2\r 0$ in the sense that we can find non-null vectors
$(\bv_s,B_s)\in \cV^{\perp_{Q_s}}$ such that $\|(\bv_s,B_s) - ({\bf 0},1)\|_{\mathds{H}}=O(s^2)$. 
Since then $\beta_s[(\bv_s,B_s)]\r 0$ and $(1/s^2)\alpha[(\bv_s, B_s)] \r 0$ as $s\r 0$, one sees from (\ref{hatQ}) 
that $(1/s^2)\hat{Q}_s[(\bv_s,B_s)]\r 2\kappa$ as $s\r 0$. Thus $Q_s[(\bv_s,B_s)]$ has the same sign as $\kappa$ for  
$s$ sufficiently small. On the other hand, for all $s\in (0,\epsilon]$ the form $Q_s$ is non-degenerate and, thus, $Q_s[(\bv_s, B_s)]\neq 0$.
Due to continuity, it follows that $Q_s[(\bv_s, B_s)]$ has the same sign as $\kappa$ for all $s\in (0,\epsilon]$.  
\qed

\subsection{The Lagrangian path $\ell(t)$ as a symplectic reduction.}

\noindent As in Sec.\ \ref{sectionsymplectic}, let $\Phi_t$ be the derivative at $(p,v)$ of the Hamiltonian flow 
$\zeta_t$. The vertical distribution $VTM$ on $TM$, given by the tangent spaces to the fibers of 
$TM\r M$, is Lagrangian. By pulling it back by $\Phi_t$, we obtain a Lagrangian path
\begin{equation}\nonumber
\check{\ell}: [0,T] \r \Lambda \bigl(T_{(p,v)}TM \bigr), \quad \check{\ell}(t) = (\Phi_t)^{-1} V_{(\gamma,\dot{\gamma})}TM.
\end{equation}

\noindent The previously defined Lagrangian path $\ell(t)$ is nothing but the symplectic reduction of $\check{\ell}(t)$ by
the coisotropic subspace $T_{(p,v)} E^{-1}(\kappa)\subset T_{(p,v)} TM$ in the sense of \cite{vitorio}. Indeed, 
$\langle S \rangle$ is the symplectic-orthogonal complement of $T E^{-1}(\kappa)$ and we have $V E^{-1}(\kappa)=VTM \cap T E^{-1}(\kappa)$, 
so that the Lagrangian subbundle $\bar{L}$ of
$TE^{-1}(\kappa)/\langle S \rangle$ is the symplectic reduction of the Lagrangian distribution $L$ along $E^{-1}(\kappa)$ given by
\begin{equation}\nonumber
L = \langle S \rangle + VTM\cap T E^{-1}(\kappa).
\end{equation} 
It is clear that $\check{\ell}(t)\cap \langle S_{(p,v)}\rangle=\{0\}$ for all $t$ since $S$ is never vertical and 
$\Phi_t(S_{(p,v)}) = S_{(\gamma,\dot{\gamma})}$. Therefore, we can apply Theorem 2.3 of \cite{vitorio} to conclude:

\begin{proposition}\label{propreduction}
The energy-constrained Maslov index $\mu_\kappa(\gamma)$ is equal to the Maslov index of $\check{\ell}|_{[\epsilon,T]}$ 
with respect to the reference Lagrangian $L_{(p,v)}$.
\end{proposition} 

\subsection{Computation of Kashiwara indices.}

\noindent Intersections of the Lagrangian path $\check{\ell}(t)$ with the Maslov cycle $\Sigma\bigl( V_{(p,v)}TM \bigr)$ correspond to 
ordinary conjugate instants along $\gamma$; this follows from the description of $\Phi_t$ in terms of ordinary Jacobi fields.
We shall call the Maslov index of $\check{\ell}|_{[\epsilon,T]}$ with respect to the reference Lagrangian 
$V_{(p,v)}TM$ the {\sl ordinary Maslov index} of $\gamma$ and denote it by $\mu(\gamma)$. 

According to Eq.\ (\ref{kashiwaraformula}),
the difference of the Maslov indices of $\check{\ell}|_{[\epsilon,T]}$ with respect to the reference Lagrangians $L_{(p,v)}$ and $V_{(p,v)}TM$, respectively,
can be computed as half the difference of Kashiwara indices. Together with Proposition \ref{propreduction}, this gives us
\begin{equation}\label{differencemaslovindices}
\mu_\kappa(\gamma) - \mu(\gamma) =   \frac{1}{2}\tau\bigl( L_{(p,v)} , V_{(p,v)}TM , \check{\ell}(T) \bigr) 
- \frac{1}{2}\tau\bigl( L_{(p,v)} , V_{(p,v)}TM, \check{\ell}(\epsilon) \bigr).
\end{equation}
For the rest of this section, let $s=\epsilon, T$. Also, we shall omit the $s$ in $\gamma(s)$ and $\dot{\gamma}(s)$.

Since the Kashiwara index is invariant by symplectic transformations, 
by applying $\Phi_s$ we obtain that $\tau \bigl( L_{(p,v)} , V_{(p,v)}TM , \check{\ell}(s)  \bigr) = \tau \bigl( \Phi_s L_{(p,v)} , \Phi_s V_{(p,v)}TM , V_{(\gamma,\dot{\gamma})}TM \bigr)$.
Let us denote the Lagrangian triplet on the right-hand side by $(L_1,L_2,L_3)$. For the following, have in mind the identification (\ref{identificationtangentspace}). 
Since $s$ is not an ordinary conjugate instant along $\gamma$, $L_2$ has null intersection with $V_{(\gamma,\dot{\gamma})}TM=\{0\}\oplus T_\gamma M$. 
Therefore there exists a linear map $A_s:T_\gamma M \r T_\gamma M$ such that $L_2 = \bigl\{(u,A_su) \; : \; u\in T_\gamma M \bigr\}$.

\begin{lemma}\label{lemmaadjoint}
The adjoint $A_s^*$ of $A_s$ with respect to $g(\gamma)$ is equal to $A_s-Y$.
\end{lemma}
\Proof
On one hand, for all $u_1,u_2\in T_\gamma M$ we have 
\begin{eqnarray*}
\omega(\gamma)\bigl[ (u_1,A_su_1), (u_2, A_su_2) \bigr] & = & g(\gamma)\bigl[ A_su_1 , u_2 \bigr] - g(\gamma)\bigl[ u_1, A_su_2 \bigr] + 
g(\gamma)\bigl[ u_1, Y[u_2] \bigr] \\
& = & g(\gamma) \bigl[ u_1 , A_s^* u_2 - A_s u_2 + Y[u_2]   \bigr].
\end{eqnarray*}
On the other hand, since $L_2$ is Lagrangian, the above must be zero.
\qed

\begin{lemma}
We have $A_s\dot{\gamma}-Y[\dot{\gamma}]\neq 0$ and the Lagrangian subspace $L_1\subset T_\gamma M \oplus T_\gamma M$ is given by 
\begin{equation}\nonumber
L_1 = \bigl\langle (\dot{\gamma},Y[\dot{\gamma}]) \bigr\rangle + \left\{ (u, A_s u) \; : \; u \in \bigl\langle A_s\dot{\gamma}-Y[\dot{\gamma}] \bigr\rangle^\perp \right\}.
\end{equation}
\end{lemma}
\Proof
Applying $\Phi_s$ to the equality $L_{(p,v)}= \langle S_{(p,v)} \rangle + \{0\}\oplus \langle v \rangle^\perp$ we obtain
$L_1=\Phi_s L_{(p,v)}=\langle S_{(\gamma,\dot{\gamma})} \rangle + \Phi_s ( \{0\} \oplus \langle v \rangle^\perp )$. We have
$S_{(\gamma,\dot{\gamma})}=(\dot{\gamma}, Y[\dot{\gamma}])$ and, since  
$\Phi_s(T_pM\oplus \langle v \rangle^\perp) = T_\gamma M \oplus \langle \dot{\gamma}\rangle^\perp$, the subspace 
$\Phi_s(\{0\}\oplus \langle v \rangle^\perp)$ is the intersection of $L_2$ with $T_\gamma M\oplus\langle \dot{\gamma} \rangle^\perp$, which in turn
is given by $\bigl\{ (u,A_su) \; : \; u\in T_\gamma M \; \mbox{and} \; g(\gamma)[A_su,\dot{\gamma}]=0 \bigr\}$. By the Lemma \ref{lemmaadjoint}, the equation
$g(\gamma)[A_su,\dot{\gamma}]=0$ is the same as $g(\gamma)[u, A_s\dot{\gamma} - Y[\dot{\gamma}]]=0$. 
If we had $A_s\dot{\gamma}-Y[\dot{\gamma}]= 0$, then $L_1$ would have dimension $n+1$, an absurd.
\qed

We have therefore an isomorphism
\begin{equation}\nonumber
\left\{
\begin{aligned}
& \R\oplus \bigl\langle A_s \dot{\gamma} - Y[\dot{\gamma}] \bigr\rangle^\perp \oplus T_\gamma M \oplus T_\gamma M \stackrel{\simeq}{\longrightarrow}  L_1\oplus L_2 \oplus L_3,  \\
& (\lambda , u_1,u_2,u_3) \mapsto  \bigl( \lambda(\dot{\gamma}, Y[\dot{\gamma}]) + x_1\, , \, x_2 \, , \, x_3 \bigr), 
\end{aligned}
\right.
\end{equation}

\noindent where $x_i=(u_i,A_su_i)$, for $i=1,2$, and $x_3=(0,u_3)$. As in Sec.\ \ref{sectionkashiwara}, let $\srQ$ be the quadratic form corresponding to the
Lagrangian triplet $(L_1,L_2,L_3)$ and denote by $\tilde{\srQ}$ the pull-back of $\srQ$ by the above
isomorphism. We shall describe $\tilde{\srQ}$. For $i=1,2,3$, let $u_i,v_i\in T_\gamma M$, with $u_1,v_1\in \langle A_s\dot{\gamma} - Y[\dot{\gamma}] \rangle^\perp$.
Using that $\omega(\gamma)$ vanishes on $L_2$, and that 
\begin{align*}
\omega(\gamma)\bigl[ (\dot{\gamma},Y[\dot{\gamma}]) , (w,A_s w) \bigr] & \;=\;  g(\gamma)\bigl[ Y[\dot{\gamma}] , w  \bigr] -g(\gamma)\bigl[\dot{\gamma} , A_s w \bigr] + g(\gamma)\bigl[\dot{\gamma}, Y[w]\bigr]  && \\
& \; = \;  - g(\gamma)\bigl[ w , A_s\dot{\gamma} - Y[\dot{\gamma}] \bigr] && \forall w\in T_\gamma M,
\end{align*}
a direct computation yields 
\begin{equation}\label{eqpullbackquadraticform}
\begin{aligned}
2\tilde{\srQ}\bigl[ (\lambda_1, u_1,u_2,u_3), (\lambda_2, v_1,v_2,v_3) \bigr]   & =     
-\lambda_1 \left( g(\gamma)\bigl[ v_2 , A_s \dot{\gamma} - Y[\dot{\gamma}] \bigr]  - g(\gamma)[v_3,\dot{\gamma}] \right) \\
 & \quad -\lambda_2 \left( g(\gamma)\bigl[ u_2 , A_s \dot{\gamma} - Y[\dot{\gamma}] \bigr]  - g(\gamma)[u_3,\dot{\gamma}] \right) \\
 & \quad -\left( g(\gamma)[u_2 , v_3] + g(\gamma)[v_2, u_3] \right) + g(\gamma)[u_3,v_1] + g(\gamma)[v_3,u_1].
\end{aligned}
\end{equation}
In order to introduce convenient coordinates to describe the matrix of $\tilde{\srQ}$, it is useful to distinguish 
between the following three cases,
\begin{equation}\label{eqcases}
\begin{aligned}
1) & \quad\; g(\gamma)\bigl[ A_s\dot{\gamma} - Y[\dot{\gamma}] \, , \,  A_s\dot{\gamma} - Y[\dot{\gamma}] \bigr]<0, \\
2) & \quad\; g(\gamma)\bigl[ A_s\dot{\gamma} - Y[\dot{\gamma}] \, , \,  A_s\dot{\gamma} - Y[\dot{\gamma}] \bigr]>0, \\
3) & \quad\; g(\gamma)\bigl[ A_s\dot{\gamma} - Y[\dot{\gamma}] \, , \, A_s\dot{\gamma} - Y[\dot{\gamma}] \bigr]=0. 
\end{aligned}
\end{equation}

\noindent Suppose, for instance, that we are in the third case above. In this case, we can introduce  
a $g-$orthonormal basis $e_1$, ..., $e_n$ on $T_\gamma M$ such that 
\begin{equation}\nonumber
A_s\dot{\gamma} - Y[\dot{\gamma}] = \lambda(e_1+e_n), \quad \mbox{for some $\lambda\neq 0$}.
\end{equation}

\noindent Recall that by $g-$orthonormal we mean that $(g(\gamma)[e_i,e_j])$ is the matrix
$I_{n,p}$ introduced in (\ref{matrixInp}).
A basis for $\bigl\langle A_s \dot{\gamma} - Y[\dot{\gamma}] \bigr\rangle^\perp$
is then $e_1+e_n$, $e_2$, ..., $e_{n-1}$. We equip the space
$\R\oplus \bigl\langle A_s \dot{\gamma} - Y[\dot{\gamma}] \bigr\rangle^\perp \oplus T_\gamma M \oplus T_\gamma M$ 
with the direct sum of the bases at hand. Also, let $a_1$, ..., $a_n$ be such that $\dot{\gamma} = \sum_{i=1}^n a_i e_i$. 
From (\ref{eqpullbackquadraticform}) it follows that the matrix of $2\tilde{\srQ}$ in that basis is

\newpage

\begin{equation}\nonumber
\AddToShipoutPictureBG*{
\begin{tikzpicture}[overlay,remember picture]
\usetikzlibrary{fit}
 \node[fill=blue!10,fit=(n1)(n2)]{};
 \node[fill=red!10,fit=(n3)(n4)]{};
 \node[fill=red!10,fit=(n5)(n6)]{};
 \node[fill=blue!10,fit=(n7)(n8)]{};
\end{tikzpicture}}
\Upsilon={\footnotesize \left(\begin{smallmatrix}
0 & & & 0 & 0 & \cdots & 0 & 0 & 0 & & & \lambda & 0 & \cdots & 0 & 0 & -\lambda & & & -a_1 & \cdots & -a_p & a_{p+1} & \cdots & a_n \\ 
&&&&&&&&&&&&&&&&&&&&&& && \\ 
&&&&&&&&&&&&&&&&&&&&& && \\
0 & & & 0 & 0 & \cdots & 0 & 0 & 0 & & &    0   & 0 & \cdots & 0 & 0 & 0 & & &   \tikzmarknode{n1}{-1}  & 0   &  0   & \cdots & 0       & 1  \\
0 & & & 0 & 0 & \cdots & 0 & 0 & 0 & & &    0   & 0 & \cdots & 0 & 0 & 0 & & &   0  &  -1   & 0   & \cdots & 0       & 0 \\
\cdots & & & \cdots & \cdots & \cdots & \cdots & \cdots & \cdots   & & &  \cdots  & \cdots   & \cdots & \cdots & \cdots & \cdots & & & \cdots  & \cdots   & \cdots & \cdots & \cdots & \cdots  \\
\cdots & & & \cdots & \cdots & \cdots & \cdots & \cdots & \cdots   & & &  \cdots  & \cdots   & \cdots & \cdots & \cdots & \cdots & & & \cdots  & \cdots   & \cdots & \cdots & \cdots & \cdots  \\
0 & & & 0 & 0 & \cdots & 0 & 0 & 0 & & &    0    & 0 & \cdots & 0 & 0 & 0 & & &   0  &  0   & 0   & \cdots & 1       & \tikzmarknode{n2}{0} \\
&&&&&&&&&&&&&&&&&&&&&& && \\
&&&&&&&&&&&&&&&&&&&&&& && \\
\lambda & & & 0 & 0 & \cdots & 0 & 0 & 0 & & &    0    & 0 & \cdots & 0 & 0 & 0 & & &   \tikzmarknode{n3}{1}  &  0   & 0   & \cdots & 0       & 0 \\
0       & & & 0 & 0 & \cdots & 0 & 0 & 0 & & &    0    & 0 & \cdots & 0 & 0 & 0 & & &   0  &  1   & 0   & \cdots & 0       & 0 \\
0       & & & 0 & 0 & \cdots & 0 & 0 & 0 & & &    0    & 0 & \cdots & 0 & 0 & 0 & & &   0  &  0   & 1   & \cdots & 0       & 0 \\
\cdots & & & \cdots & \cdots & \cdots & \cdots & \cdots & \cdots  & & &    \cdots    & \cdots & \cdots & \cdots & \cdots & \cdots & & & \cdots  & \cdots   & \cdots & \cdots & \cdots & \cdots \\
0       & & & 0 & 0 & \cdots & 0 & 0 & 0  & & &    0    & 0 & \cdots & 0 & 0 & 0 & & &   0  &  0   & 0   & \cdots & -1       & 0 \\
-\lambda       & & & 0 & 0 & \cdots & 0 & 0 & 0 & & &    0    & 0 & \cdots & 0 & 0 & 0 & & &   0  &  0   & 0   & \cdots &  0       & \tikzmarknode{n4}{-1} \\  
&&&&&&&&&&&&&&&&&&&&&& && \\
&&&&&&&&&&&&&&&&&&&&&& && \\
-a_1    & & & \tikzmarknode{n7}{-1} & 0 & \cdots & 0 & 0 & 0 & & &    \tikzmarknode{n5}{1}    & 0 & \cdots & 0 & 0 & 0 & & &   0  &  0   & 0   & \cdots &  0       & 0 \\ 
\cdots  & & & 0 & -1 & \cdots & 0 & 0 & 0                    & & &    0  & 1 & \cdots & 0 & 0  & 0 & & &   0  &  0   & 0   & \cdots &  0       & 0 \\
-a_p    & & & \cdots & \cdots & \cdots & \cdots & \cdots & \cdots & & &   \cdots & \cdots & \cdots & \cdots & \cdots & \cdots & & &   \cdots  &  \cdots   & \cdots   & \cdots &  \cdots       & \cdots \\
a_{p+1} & & & 0 &0 & \cdots & 1 & 0 & 0 & & &  0 & 0 & \cdots & -1 & 0 & 0 & & &   0  &  0   & 0   & \cdots &  0       & 0 \\
\cdots & & & 0 &0 & \cdots & 0 & 1 & 0 & & &  0 & 0 & \cdots &   0 &-1 & 0 & & &   0  &  0   & 0   & \cdots &  0       & 0 \\
a_{n}   & & & 1 &0 & \cdots & 0 & 0 & \tikzmarknode{n8}{0} & & &    0 & 0 & \cdots &  0 & 0 & \tikzmarknode{n6}{-1} & & &   0  &  0   & 0   & \cdots &  0       & 0 \\
\end{smallmatrix}\right),}
\end{equation}

\noindent where:
\begin{itemize}
\item[--] the red blocks are equal to $-I_{n,p}$,
\item[--] the upper blue block is $(n-1)\times n$ and has a diagonal formed by $p$ minus ones followed by $n-p-1$ plus ones,
\item[--] the bottom blue block is the transpose of the upper one.
\end{itemize}
By performing a few elementary cogredient transformations (i.e., a sequence of row operations followed by the same operations on the columns) on 
$\Upsilon$ we arrive at the matrix
\begin{equation}\nonumber
\Upsilon' = 
\begin{pmatrix}
2\lambda(a_1-a_n) &  &  &   \\
                  & {\rm O}_{(n-1)\times (n-1)} &   &    \\                  
                  &               & {\rm O}_{n\times n} & - I_{n,p} \\                  
                  &               & -I_{n,p}  & {\rm O}_{n\times n}  
\end{pmatrix},
\end{equation} 
where ${\rm O}_{k\times k}$ is a $k\times k$ matrix of zeros and the omitted elements are zeros. The $2n\times 2n$ 
bottom-right block of $\Upsilon'$ has null signature. Also, it follows from $A_s\dot{\gamma}-Y[\dot{\gamma}]=\lambda(e_1+e_n)$ and Lemma \ref{lemmaadjoint}
that $2\lambda(a_1-a_n)=-2g(\gamma)[A_s\dot{\gamma} , \dot{\gamma}]$. Therefore,
\begin{equation}\nonumber
{\rm sig}\, \srQ = {\rm sig}\,\tilde{\srQ} = {\rm sig}\, \Upsilon' = 
\left\{
\begin{aligned}
1, & \quad {\rm if} \quad g(\gamma)[A_s\dot{\gamma} , \dot{\gamma}]< 0, \\
-1,&\quad {\rm if} \quad g(\gamma)[A_s\dot{\gamma} , \dot{\gamma}]> 0, \\
0, & \quad {\rm if} \quad g(\gamma)[A_s\dot{\gamma} , \dot{\gamma}]= 0.
\end{aligned}
\right.
\end{equation}
The cases 1) and 2) in (\ref{eqcases}) can be treated quite similarly and provide the same value for ${\rm sig}\, \srQ$ as above.
Let now $\check{J}_s^*$ be as in Proposition \ref{propindexperp}. Since $L_2=\Phi_s V_{(p,v)}TM$ consists of all pairs $(\check{J},D\check{J}/dt)|_{t=s}$ 
where $\check{J}$ is an ordinary Jacobi field vanishing at $t=0$, then $A_s\dot{\gamma} = D\check{J}_s^*/dt|_{t=s}$. In particular, 
$g(\gamma)[A_s\dot{\gamma},\dot{\gamma}]\neq 0$. 
In conclusion, we have 
\begin{proposition}\label{propkashiwaracomputation}
Let $s=\epsilon, T$, and let $\check{J}_s^*$ be as in Proposition \ref{propindexperp}. Then,
\begin{equation}\nonumber
\tau\bigl( L_{(p,v)} , V_{(p,v)}TM , \check{\ell}(s) \bigr) = 
\left\{
\begin{aligned}
1, & \quad {\rm if} \quad g(\gamma)\left[ D \check{J}^*_s/dt , \dot{\gamma} \right]\big|_{t=s}<0, \\
-1, & \quad {\rm if} \quad g(\gamma)\left[ D \check{J}^*_s/dt , \dot{\gamma} \right]\big|_{t=s}>0.
\end{aligned}
\right.
\end{equation}
\end{proposition}

\begin{remark}
It follows from (\ref{differencemaslovindices}) and the above computations that, for a given instant $t_0$ for which there are no energy-constrained or
ordinary conjugate instants in a neighborhood of $t_0$ (except possibly for $t_0$), then 
$\big|\mu_\kappa(\gamma, t_0) - \mu(\gamma, t_0) \big| \leq 1$ (recall the definition of Maslov index across an instant).
Indeed, in the notation of Proposition \ref{propindexperp}, for $\delta>0$ small enough we have
\begin{equation}\nonumber
\mu_\kappa(\gamma,t_0) - \mu(\gamma,t_0) = {\rm ind}\, Q_{t_0+\delta}\big|_{\cV^{\perp_{Q_{t_0+\delta}}}} - 
{\rm ind}\, Q_{t_0-\delta}\big|_{\cV^{\perp_{Q_{t_0-\delta}}}},
\end{equation} 
which is 0 or $\pm 1$ since $\cV$ is one-codimensional. In the case $t_0$ is non-degenerate as a conjugate instant (energy-constrained and/or
ordinary), the difference $\mu_\kappa(\gamma,t_0) - \mu(\gamma,t_0)$ is just the difference of the respective signatures of $t_0$.  
\end{remark}

\subsection{Proof of Theorem \ref{morseindextheorem}.}

\noindent By comparing (\ref{eqdifferencespectralflow}) and Proposition \ref{propindexperp} with
(\ref{differencemaslovindices}) and Proposition \ref{propkashiwaracomputation}, we obtain
\begin{equation}\label{eq7834}
{\rm sf}_\kappa(\gamma) - {\rm sf}(\gamma) = \mu(\gamma) - \mu_\kappa(\gamma).
\end{equation}
The proof of Theorem \ref{morseindextheorem} is thus reduced to the proof of the equality
\begin{equation}\nonumber
{\rm sf}(\gamma) = - \mu(\gamma),
\end{equation}
which in turn is the usual Morse index theorem proved in Proposition 2.5 of \cite{portaluri_piccione}.

\subsection{Proof of Theorem \ref{theorembifurcationmaslov} and Corollary \ref{corollarymaslovacross}.}\label{subsectionproofbifurcation}
\noindent Theorem \ref{theorembifurcationmaslov} will follow from Theorem \ref{thmabstractbifurcation} and 
Theorem \ref{morseindextheorem} once we frame our bifurcation problem within the abstract setting of Sec.\ \ref{sectionabstractbifurcation}.
In order to do so, recall that, for each $s\in (0,T]$, $x_s : [0,1] \r M$ is the curve $x_s(t)=\gamma(st)$. We have thus    
a one-parameter family of local magnetic $\kappa-$action functionals 
centered at $x_s$,  
\begin{equation}\nonumber
\cA_{\kappa,\sigma,x_s} : \cO_s\times \R_+ \subset \Omega_{p,\gamma(s)}([0,1]) \times \R_+ \r \R, \quad s\in (0,T],
\end{equation}
each of which possessing $(x_s,s)$ as a critical point.
By means of the parametrizations (\ref{parametrization}) and the fixed frame $E_1$, ..., $E_n$ along $\gamma$, we can identify 
$\{\cA_{\kappa,\sigma,x_s}\}$ with
a smooth family of functionals with $s-$independent domains, 
\begin{equation}\nonumber
\cF_s : H^1_0([0,1],B^n_\rho)\times \R_+ \r \R,  \quad s\in  (0,T],
\end{equation}   
all of which having $({\bf 0},1)$ as the critical point corresponding to $(x_s,s)$. More precisely: 
for each $s\in(0,T]$, we use the frame $E_1^s$, ..., $E_n^s$ along $x_s$ given by $E_i^s(t)=E_i(st)$
to define the map $\psi_s: [0,1]\times B^n_\rho \r M$ as in (\ref{mapparametrization}) with $x_*$ replaced by $x_s$.
Observe that $\psi_s(t,\bx)=\psi_T((s/T)t,\bx)$ and, hence, $(s,t,\bx)\mapsto \psi_s(t,\bx)$ is smooth. If we choose $\rho$ so that
$\psi_T(t, \cdot):B_\rho^n \r M$ is a diffeomorphism onto an open set for each $t$, the same will be true of the maps 
$\psi_s(t,\cdot):B_\rho^n \r M$ for all $s$ and $t$. Let then $\Psi_s:H^1_0([0,1],B^n_\rho) \r \Omega_{p,\gamma(s)}([0,1])$ be 
the parametrization around $x_s$ corresponding to $\psi_s$ as in (\ref{parametrization}). 
The desired family $\{\cF_s\}_{s\in (0,T]}$ is then defined by 
\begin{equation}\nonumber
\cF_s (\bx , r) = \cA_{\kappa,\sigma,x_s}(\Psi_s(\bx) , sr). 
\end{equation}
For the smoothness of $(s,\bx,r)\mapsto \cF_s(\bx,r)$, let us split $\cA_{\kappa,\sigma,x_s}$ as $\cA_\kappa^s + \cC_{\sigma,x_s}$, for 
$\cA_\kappa^s:\Omega_{p,\gamma(s)}([0,1])\times \R_+ \r \R$ the 
$\kappa-$action functional of the semi-Riemannian Lagrangian. The smoothness of $(s,\bx,r)\mapsto \cA_\kappa^s(\Psi_s(\bx),sr)$
follows easily from that of $(s,t,\bx)\mapsto \psi_s(t,\bx)$. Furthermore, 

\begin{lemma}
The map $(s,\bx)\mapsto \cC_{\kappa,x_s}(\Psi_s(\bx))$ is smooth.
\end{lemma}
\Proof
As in the discussion following Corollary \ref{corollarycriticalpoint}, let $\theta$ be a $1-$form on $[0,1]\times B_\rho^n$
such that $\psi_T^*\sigma=d\theta$. For each $s$, let $h_s:[0,1]\times B_\rho^n \r [0,1]\times B_\rho^n$ be the diffeomorphism
$h_s(t,\bx)=((s/T)t,\bx)$. Since $\psi_s = \psi_T\circ h_s$, we have that 
$\psi_s^*\sigma=d(h_s^*\theta)$. Thus, according to Eq. (\ref{eqlocalexpression}),
\begin{equation}\nonumber
\cC_{\sigma,x_s}(\Psi_s(\bx)) = \int_0^1(h_s^*\theta)(t,\bx)[(1,\dot{\bx})] dt - \int_0^1 (h_s^*\theta)(t,{\bf 0})[(1,{\bf 0})]dt.
\end{equation}
The claimed smoothness follows now from the smoothness of $(s,t,\bx)\mapsto (h_s^*\theta)(t,\bx)$. 
\qed

Each functional $\cF_s$ has $({\bf 0},1)$ as a critical point since $\Psi_s({\bf 0})=x_s$.
Also, since the derivative of the diffeomorphism $({\bf x},r)\mapsto(\Psi_s({\bf x}),sr)$ at $({\bf 0},1)$ is the isomorphism (\ref{isomorphismtrivialization}), the second derivative $(D^2\cF_s)_{({\bf 0},1)}$ is equal to the quadratic form $Q_s$ of 
Sec. \ref{sectionanalytic}. On the other hand, it is clear that
an instant $t_0$ is an energy-constrained bifurcation instant along $\gamma$ if, and only if, it is a bifurcation instant
for the family $\{\cF_s\}_{s\in(0,T]}$ in the sense of Sec.\ \ref{sectionabstractbifurcation}; indeed, observe that
if $s_n$ and $(\bx_n,r_n)$ are sequences as in the definition in Sec.\ \ref{sectionabstractbifurcation}, then
$t_n=s_n$ and $t_n'=r_n s_n$ are sequences as in Definition \ref{definitionbifurcation}.
In particular, such an
instant must be a degenerate instant for $\{Q_s\}$ and thus an energy-constrained conjugate instant along $\gamma$.
Theorem \ref{theorembifurcationmaslov} follows now from Theorem \ref{thmabstractbifurcation} and 
Theorem \ref{morseindextheorem}. 

As for the proof of Corollary \ref{corollarymaslovacross}, observe that Theorem \ref{theorembifurcationmaslov} together with
the additivity under concatenation of both the Maslov index and the spectral flow imply the equality 
${\rm sf}(Q,t_0)=-\mu_\kappa(\gamma,t_0)$. We now apply Corollary \ref{corollaryabstractbifurcationacross} to obtain
the result.
 
\subsection{Computing the difference ${\rm sf}_\kappa(\gamma)-{\rm sf}(\gamma)$ in terms of an {\sl orbit strip}. }\label{subsectionorbitcylinder}

\noindent The concept of {\sl orbit cylinder} has played an important role in the computation of
the difference ${\rm sf}_\kappa(\gamma)-{\rm sf}(\gamma)$ in the works \cite{paternain_merry}, \cite{paternain2015} and 
\cite{portaluri_wu_yang}. In the present fixed endpoints case, an analogous concept would be the following:  

\begin{definition}
Let $\gamma:[0,T]\r M$ be an electromagnetic geodesic with energy $\kappa\neq 0$ connecting points $p$ and $q$. An 
{\sl orbit strip} around it is a smooth variation of $(x,T)$ in $\Omega_{p,q}([0,1])$, $(-\epsilon, \epsilon)\ni s \mapsto (x_s,T_s)$, 
such that $\gamma_s(t)=x_s(t/T_s)$ is an electromagnetic geodesic with energy $E(\gamma_s,\dot{\gamma}_s)\equiv \kappa+s$ for all $s$. 
\end{definition}

For the next proposition, let $T_0'$ denote $dT_s/ds|_{s=0}$. Also, for a non null real number $A$, we let 
${\rm sign}\, A$ be $\pm1$ according to $A>0$ or $A<0$, respectively. 

\begin{proposition}
Let $T$ not be an ordinary conjugate instant along $\gamma$. If there exists an orbit strip 
$(-\epsilon, \epsilon)\ni s \mapsto (x_s,T_s)$ around $\gamma$, then $T$ is not an energy-constrained conjugate 
instant along $\gamma$, $T_0'\neq 0$ and
\begin{equation}\nonumber
{\rm sf}_\kappa(\gamma) - {\rm sf}(\gamma) = \frac{1}{2}{\rm sign}(-T_0') - \frac{1}{2}{\rm sign}\, \kappa .
\end{equation}
\end{proposition}
\Proof
Let $\check{J}$ be the ordinary Jacobi field along $\gamma$ given by $\check{J}=d\gamma_s/ds|_{s=0}$. From the energy assumption, we obtain
\[
1\equiv \frac{d}{ds}\Big|_{s=0} \frac{1}{2}g(\gamma)\left[ \dot{\gamma}_s, \dot{\gamma}_s\right] \equiv g(\gamma)\left[ D\check{J}/dt , \dot{\gamma} \right].
\]
In particular, $\check{J}\neq 0$. On the other hand, by a simple computation, 
$\check{J}(t)=V(t/T) - (T_0'/T)t \dot{\gamma}(t)$, where $V=dx_s/ds|_{s=0}$. Since $V(0)=0$
and $V(1)=0$, we have $\check{J}(0)=0$ and $\check{J}(T)=-T_0'\dot{\gamma}(T)$. 
Therefore, $T_0'\neq 0$, since $T$ is not an ordinary conjugate instant, and $-(1/T_0')\check{J}$
is equal to the ordinary Jacobi field $\check{J}^*_T$ of Proposition \ref{propindexperp}. Since then
$g\left[ D\check{J}^*_T/dt , \dot{\gamma} \right] = -(1/T_0')g\left[ D\check{J}/dt , \dot{\gamma} \right]=-1/T_0'$, 
it follows from the same proposition that 
\[
{\rm ind}\, Q_T\big|_{\cV^{\perp_{Q_T}}} = \frac{1}{2} - \frac{1}{2}{\rm sign}(-T_0')
\]
and that $T$ is not an energy-constrained conjugate instant. We can thus apply formula (\ref{eqdifferencespectralflow})
and conclude the result from the above equality and from Proposition \ref{propepsilonindex}. 
\qed

\bigskip 
\medskip

\centerline{\headfont Acknowledgements}
\medskip

I would like to thank the anonymous Referee for all the constructive suggestions that helped improve the quality and clarity of the work, as well for 
pointing out the shortcomings of a previous version of the same.

\bigskip
\centerline{\headfont Data availability statement}
\medskip

No new data were created or analyzed in this study. Data sharing is not applicable to this article.

\bigskip
\centerline{\headfont References}
\medskip
\bibliography{vitoriobifurcation_v3}

\begin{thebibliography}{10}

\bibitem{abbondandolo}
A.~Abbondandolo.
\newblock Lectures on the free period {L}agrangian action functional.
\newblock {\em J. Fixed Point Theory Appl.}, 13(2):397--430, 2013.

\bibitem{paternain2015}
A.~Abbondandolo, L.~Macarini, and G.~P. Paternain.
\newblock On the existence of three closed magnetic geodesics for subcritical
  energies.
\newblock {\em Comment. Math. Helv.}, 90(1):155--193, 2015.

\bibitem{gabrielle}
L.~Asselle and G.~Benedetti.
\newblock The {L}usternik-{F}et theorem for autonomous {T}onelli {H}amiltonian
  systems on twisted cotangent bundles.
\newblock {\em J. Topol. Anal.}, 8(3):545--570, 2016.

\bibitem{assenza}
V.~Assenza.
\newblock Magnetic curvature and existence of a closed magnetic geodesics on
  low energy levels.
\newblock {\em Int. Math. Res. Not. IMRN}, (21):13586--13610, 2024.

\bibitem{assenza2}
V.~Assenza, J.~M. Reber, and I.~Terek.
\newblock Magnetic flatness and {E}. {H}opf's theorem for magnetic systems.
\newblock {\em Comm. Math. Phys.}, 406(2), 2025.

\bibitem{patodi}
M.~Atiyah, V.~Patodi, and I.~Singer.
\newblock Spectral asymmetry and {Riemannian} geometry. {III}.
\newblock {\em Math. Proc. Camb. Philos. Soc.}, 79:71--99, 1976.

\bibitem{benevieri}
P.~Benevieri and P.~Piccione.
\newblock On a formula for the spectral flow and its applications.
\newblock {\em Math. Nachr.}, 283(5):659--685, 2010.

\bibitem{piccione_notices}
R.~Bettiol and P.~Piccione.
\newblock Instability and bifurcation.
\newblock {\em Notices Amer. Math. Soc.}, 67(11):1679--1691, 2020.

\bibitem{cappellmillerlee}
S.~E. Cappell, R.~Lee, and E.~Y. Miller.
\newblock On the {Maslov} index.
\newblock {\em Commun. Pure Appl. Math.}, 47(2):121--186, 1994.

\bibitem{spectralflowbook}
N.~Doll, H.~Schulz-Baldes, and N.~Waterstraat.
\newblock {\em Spectral flow---a functional analytic and index-theoretic
  approach}.
\newblock De Gruyter Studies in Mathematics. De Gruyter, Berlin, 2023.

\bibitem{recht}
P.~M. Fitzpatrick, J.~Pejsachowicz, and L.~Recht.
\newblock Spectral flow and bifurcation of critical points of
  strongly-indefinite functionals. {I}. {G}eneral theory.
\newblock {\em J. Funct. Anal.}, 162(1):52--95, 1999.

\bibitem{foulon}
P.~Foulon.
\newblock G\'{e}om\'{e}trie des \'{e}quations diff\'{e}rentielles du second
  ordre.
\newblock {\em Ann. Inst. H. Poincar\'{e} Phys. Th\'{e}or.}, 45(1):1--28, 1986.

\bibitem{piccione_nonlinearity}
R.~Giamb\`o, F.~Giannoni, and P.~Piccione.
\newblock Gravitational lensing in general relativity via bifurcation theory.
\newblock {\em Nonlinearity}, 17(1):117--132, 2004.

\bibitem{giambo}
R.~Giamb\`o and M.~A. Javaloyes.
\newblock A second-order variational principle for the {L}orentz force
  equation: conjugacy and bifurcation.
\newblock {\em Proc. Roy. Soc. Edinburgh Sect. A}, 137(5):923--936, 2007.

\bibitem{partialsignatures}
R.~Giamb\`o, P.~Piccione, and A.~Portaluri.
\newblock Computation of the {M}aslov index and the spectral flow via partial
  signatures.
\newblock {\em C. R. Math. Acad. Sci. Paris}, 338(5):397--402, 2004.

\bibitem{jost}
J.~Jost and X.~Li-Jost.
\newblock {\em Calculus of Variations}, volume~64 of {\em Cambridge Studies in
  Advanced Mathematics}.
\newblock Cambridge University Press, Cambridge, 1998.

\bibitem{vergne}
G.~Lions and M.~Vergne.
\newblock {\em The Weil Representation, Maslov Index and Theta Series}.
\newblock Progr. Math. Birkh{\"a}user, 1980.

\bibitem{paternain_merry}
W.~J. Merry and G.~P. Paternain.
\newblock Index computations in {R}abinowitz {F}loer homology.
\newblock {\em J. Fixed Point Theory Appl.}, 10(1):87--111, 2011.

\bibitem{portaluri2007}
M.~Musso, J.~Pejsachowicz, and A.~Portaluri.
\newblock Morse index and bifurcation of {$p$}-geodesics on semi {R}iemannian
  manifolds.
\newblock {\em ESAIM Control Optim. Calc. Var.}, 13(3):598--621, 2007.

\bibitem{phillips}
J.~Phillips.
\newblock Self-adjoint {Fredholm} operators and spectral flow.
\newblock {\em Can. Math. Bull.}, 39(4):460--467, 1996.

\bibitem{portaluri_piccione}
P.~Piccione and A.~Portaluri.
\newblock A bifurcation result for semi-{R}iemannian trajectories of the
  {L}orentz force equation.
\newblock {\em J. Differential Equations}, 210(2):233--262, 2005.

\bibitem{tausk2004}
P.~Piccione, A.~Portaluri, and D.~V. Tausk.
\newblock Spectral flow, {M}aslov index and bifurcation of semi-{R}iemannian
  geodesics.
\newblock {\em Ann. Global Anal. Geom.}, 25(2):121--149, 2004.

\bibitem{maslovbook}
P.~Piccione and D.~V. Tausk.
\newblock {\em A student's guide to symplectic spaces, Grassmannians and Maslov
  index}.
\newblock Publ. Mat. IMPA. Rio de Janeiro: Instituto Nacional de Matem{\'a}tica
  Pura e Aplicada (IMPA), 2008.

\bibitem{portaluri_wu_yang}
A.~Portaluri, L.~Wu, and R.~Yang.
\newblock Linear instability of periodic orbits of free period {L}agrangian
  systems.
\newblock {\em Electron. Res. Arch.}, 30(8):833--2859, 2022.

\bibitem{reber_shen}
J.M. Reber and Y.~Shen.
\newblock Anosov magnetic flows on surfaces.
\newblock {\em arXiv:2406.18735}, 2024.

\bibitem{vitorio}
H.~Vit{\'o}rio.
\newblock On the {Maslov} index in a symplectic reduction and applications.
\newblock {\em Proc. Am. Math. Soc.}, 148(8):3517--3526, 2020.

\end{thebibliography}
\bibliographystyle{plain}

\bigskip

{\small Henrique Vit\'orio, Departamento de Matem\'atica, Universidade Federal de Pernambuco, Rua Jornalista An\'ibal Fernandes,
Cidade Universit\'aria, 50740560,
Recife, PE, Brasil. E-mail address:} {\smalltt henrique.vitori@ufpe.br}

\end{document}